\documentclass[review,10pt,3p]{elsarticle}
\usepackage{bm}
\usepackage{color}
\usepackage{amsfonts}
\usepackage{amsmath}
\usepackage{amssymb}
\usepackage{amsthm}
\usepackage{graphicx}
\usepackage{subfigure}
\usepackage{listings} 
\usepackage{multirow}
\usepackage{enumitem}
\usepackage{hyperref}
\newcommand{\blue}[1]{{#1}}

\newtheorem{thm}{Theorem}[section]

\newtheorem{theorem}[thm]{Theorem}
\newtheorem{lemma}[thm]{Lemma}
\newtheorem{remark}[thm]{Remark}
\newtheorem{example}[thm]{Example}
\def \blue#1{{#1}}
\definecolor{mygreen}{RGB}{0,120,80}

\newcommand{\lap}{\mathrm{\Delta}}
\newcommand{\del}[1]{}
\newcommand{\ins}[1]{#1}
\newcommand{\opecolor}[1]{}
\newcommand{\linkref}[2][]{#1~#2}
\newenvironment{auimage}[1]{}
  {}
\makeatletter
\@ifundefined{coi}{\newenvironment{coi}{}{}}{}
\@ifundefined{ack}{\newenvironment{ack}{}{}}{}
\@ifundefined{appgroup}{\newenvironment{appgroup}{}{}}{}
\@ifundefined{appsection}{\newcommand{\appsection}[1]{\section*{#1}}}{}
\makeatother
\numberwithin{equation}{section}
\usepackage{float}

\journal{***}
\begin{document}

\begin{frontmatter}



\title{Long-time stability analysis of an explicit exponential Runge--Kutta scheme for Cahn--Hilliard equations}





\author[gdut]{Jing Guo\corref{cor1}}
\ead{jingguo@gdut.edu.cn}
\cortext[cor1]{Corresponding author.}
\affiliation[gdut]{
  organization={School of Mathematics and Statistics, Guangdong University of Technology},
  addressline={Guangdong},
  city={Guangzhou},
  postcode={510006},
  country={China}
}

\begin{abstract}
In this paper, we present a rigorous long-time stability analysis of a second-order explicit exponential Runge--Kutta (ERK2) method for the Cahn--Hilliard  equation.
By employing Fourier spectral collocation in space and a two-stage ERK2 scheme in time, we construct a fully discrete numerical method and establish an energy dissipation law for the original energy.
{ The numerical solution is proven to be uniformly bounded in time in the discrete $H^1$ and $H^2$ norms, provided that the time step size is sufficiently small. An $\ell^\infty$ bound is then derived through a discrete Sobolev inequality. These bounds remove the typical a priori maximum-norm assumption required in previous energy-stability analyses and allow the energy dissipation criterion to be closed for the fully discrete scheme.}
Building on this uniform boundedness, we derive an optimal-order error estimate in the $\ell^2$ norm. The analytical framework developed here is general and can be extended to higher-order exponential integrators for a broader class of phase-field models.
\end{abstract}

\end{frontmatter}

\section{Introduction}\label{Xsec1-1}
The Cahn--Hilliard (CH) equation is one of the fundamental equations in the study of phase separation and coarsening in inhomogeneous systems, including glasses, alloys, and polymer blends at fixed temperatures, and was originally proposed by Cahn and Hilliard in \cite{cahn1958free}. In this work, we consider the CH equation in the form
\begin{equation}
    \label{eqn:ch_eq}
    \left\{
    \begin{aligned}
    &u_t = \Delta (- \varepsilon^2 \Delta u +  f(u)),\quad \bm{x} \in \Omega,\ t \in (0,\infty),\\
    &u(\bm{x},0) = u_0(\bm{x}),
    \end{aligned}
    \right.
\end{equation}
where $u(\bm{x},t) : \Omega \subset \mathbb{R}^d \rightarrow \mathbb{R}$ denotes the order parameter and represents the concentration difference between two components, and $\varepsilon$ is the diffusion coefficient related to the interface  width from physical aspects. For brevity, we impose a periodic boundary condition on $\Omega$, which is amenable to a particularly transparent numerical analysis via the Fourier spectral method. An extension to the case of a homogeneous Neumann boundary condition is straightforward. Under the aforementioned boundary conditions, the CH equation can be regarded as the $H^{-1}$ gradient flow associated with the Ginzburg--Landau energy functional:
\begin{equation}\label{eqn:ch_energy}
  E(u) = \int_\Omega \frac{\varepsilon^2}{2} |\nabla u|^2 + F(u) \mathrm{d}\bm{x},
\end{equation}
where $F(u) = \frac{1}{4}(u^2 -1)^2$ is the double-w\blue{e}ll potential and we thus derive $f(u) := F'(u) = u^3 - u$. The solution of \eqref{eqn:ch_eq} satisfies the energy dissipation law:
\begin{equation}
    \frac{\mathrm{d}}{\mathrm{d}t} E(u) = \left(\frac{\delta E(u)}{\delta u}, \frac{\mathrm{d} u}{\mathrm{d} t}\right) = - \int_\Omega |\nabla (-\varepsilon^2 \Delta u + f(u))|^2 \mathrm{d} \bm{x} \leq 0,
\label{Xeqn3-1.3}
\end{equation}
which indicates that the energy $E(u)$ is non-increasing over time.

In the existing literature, numerous numerical methods have been developed for the CH equation \cite{eyre1998unconditionally-ch,wise2010unconditionally,GuoWanWise24,GuoWanWise26,LiLiuWan24,SunZhanQian,LiQuanTang,LiQuan,furihata2001stable,ZhanLiu24}. Among these studies, considerable efforts have been devoted to constructing numerical schemes that are not only computationally efficient but also capable of preserving the energy stability of the CH equation.
By decomposing the energy functional into the difference of two convex parts, Eyre \cite{eyre1998unconditionally-ch} proposed the convex splitting scheme in 1998 and established its  energy stability.
In 2001, Furihata \cite{furihata2001stable} developed a conservative Crank--Nicolson scheme to solve the one-dimensional CH equation, and proved that this scheme inherits the key properties of the equation, including the mass conservation and the energy dissipation law, under a time step constraint.
Additionally, Wise et al. \cite{wise2010unconditionally} extended the convex splitting scheme to investigate the Cahn--Hilliard--Hele--Shaw system, and conducted a rigorous analysis of both the energy stability and convergence of their proposed scheme.
Chen et al. \cite{chen2024second} proposed a second-order {accurate} finite difference scheme for the Cahn--Hilliard--Navier--Stokes system with a logarithmic Flory--Huggins energy potential, leveraging a modified Crank--Nicolson scheme and explicit second-order Adams--Bashforth extrapolation. They further analyzed the unique solvability and positivity-preserving property of the fully discrete scheme.
More recently, Cheng et al. \cite{cheng2022third} developed a third-order backward differentiation formula (BDF)-type scheme for the CH equation by incorporating the convex splitting idea. They proved the scheme is energy stable with respect to a modified energy and demonstrated that it guarantees the uniform boundedness of the original energy.
While the convex splitting scheme and modified Crank--Nicolson scheme enable the efficient construction of  energy-stable schemes, these schemes are typically nonlinear. Consequently, they incur high computational costs, as nonlinear solvers must be employed at each time step.

To avoid nonlinear iterations while preserving the energy stability of the CH equation, implicit-explicit (IMEX) schemes, where the linear diffusion term is treated implicitly and the nonlinear reaction term is treated explicitly, have garnered increasing attention over the past decades.
Feng and Prohl \cite{feng2003numerical} established the error analysis for both a semi-discrete and a fully discrete finite element method applied to the CH equation. For the CH equation with dynamic boundary conditions dependent on the reaction rate, Meng et al. \cite{xiangjun2025second} proposed a second-order stabilized Crank-Nicolson scheme and provided the corresponding theoretical proofs of its stability and convergence.
Qiao et al. \cite{Li2023Stabilization} formulated an  energy-stable, second-order scheme for the nonlocal CH equation by coupling the BDF method with a second-order stabilization term. They further derived high-order consistency estimates, which yielded an $\ell^\infty$ bound for the numerical solution.
Additionally, Fu et al. \cite{fu2022energy} conducted an energy stability analysis for the exponential time differencing Runge--Kutta (ETDRK) scheme in the context of general phase-field models satisfying the Lipschitz condition. For classical implicit-explicit Runge--Kutta (IMEX RK) schemes, Fu et al. \cite{fu2022unconditionally} also presented a general proof showing that such schemes can  preserve the energy dissipation law when combined with stabilization techniques and under specific RK coefficient constraints. The energy stability analyses mentioned above rely either on a global Lipschitz condition imposed on the nonlinearity or {on} an $\ell^\infty$ bound assumption for the numerical solution.
To resolve these constraints in the context of thin film equations and the CH equation, Li et al. \cite{li2017second,li2017stabilization} employed harmonic analysis techniques in borderline spaces for first- and second-order linearized BDF schemes. This enabled them to establish long-time stability in the maximum norm, where the result depends on both the diffusion coefficient and the initial condition.
{In \cite{LiQiao24}, an $\ell^\infty$ bound on the numerical solution was established for a second-order scheme for a nonlocal CH equation under the condition $\tau = O(h)$, where $\tau$ and $h$ denote the temporal and spatial step sizes, respectively. Gottlieb {\rm et al.}~\cite{gottlieb2012long} analyzed the long-time stability of a semi-implicit scheme for the 2D Navier--Stokes equations, with sufficiently small $\tau$ bounded by a constant depending on the viscosity coefficient.}
{In addition, global-in-time energy stability of the phase field crystal equation in the \( H^{s} \) (\( s = 1,2 \)) norm has been established via energy-based techniques, including the study of a second-order ETDRK scheme in~\cite{li2025global} and a second-order IMEX Runge--Kutta scheme in~\cite{zhang2024second}.}
Nevertheless, a comprehensive characterization of the uniform-in-time stability in the maximum norm for energy-stable exponential Runge--Kutta schemes has not yet been studied.

The objective of this work is to construct a two-stage, second-order accurate exponential Runge--Kutta (ERK2) scheme for the CH equation, and to analyze the long-time numerical stability in $H^s$ $(s = 1,2)$ norms and $\ell^\infty$ norm. The main contributions are as follows:
\begin{itemize}
    \item A fully discrete ERK2 scheme is developed using Fourier spectral collocation, which preserves the original energy dissipation.
    \item Uniform-in-time boundedness of numerical solutions is established in discrete $H^1$, $H^2$, and $\ell^{\infty}$ norms  {when the time step size is sufficiently small},
    removing the boundedness assumption typically required in energy stability analyses.
    \item {The} energy stability is proven with a stabilization parameter that depends only on the $\ell^{\infty}$ norm of the solution, accompanied by {an} optimal-order convergence estimate.
\end{itemize}

The remainder of this paper is organized as follows. \linkref[Section]{\ref{sec:erk}} provides a brief review of the Fourier spectral collocation method and introduces the ERK2 scheme, along with preliminary estimates and an energy stability result.  In \linkref[Section]{\ref{sec:stab}}, we establish the uniform-in-time boundedness of the numerical solution in the discrete $H^1$, $H^2$ and $\ell^\infty$ norm, ensuring the long-time stability and global-in-time energy stability of the ERK2 scheme.  \linkref[Section]{\ref{sec:error}} presents the convergence analysis, where an optimal-order error estimate is  obtained.  {In \linkref[Section]{\ref{sec:numex}},  two numerical examples are implemented to validate the theoretical results regarding energy dissipation, long-time stability and convergence rate. } Concluding remarks are given in \linkref[Section]{\ref{sec:con}}.

\section{Fully discrete ERK scheme with Fourier collocation discretization}\label{Xsec2-2}\label{sec:erk}
\subsection{Review of the Fourier collocation method}\label{Xsec3-2.1}
The Fourier collocation method, also called the Fourier pseudo-spectral method, is closely related to the Fourier spectral method. Its main difference is that it approximates the Fourier coefficients using the composite trapezoidal rule. When combined with fast algorithms such as the fast Fourier transform (FFT), this approach can greatly accelerate computations.
To simplify the notation used in our analysis, we consider the three-dimensional case in this work, take the domain as $\Omega = (0,L)^3$, let $N \in \mathbb{N}^+$ be a positive even integer and define the following index sets:
\begin{align}
\label{set_SN}
   &\mathcal{S}_N = \Big\{\bm m := (p,q,r) \in \mathbb{Z}^3 | 0\leq \bm m \leq N-1 \Big\},\\
   \label{set_hatSN}
  &\widehat{\mathcal{S}}_N = \Big\{\bm k := (k, \ell, m) \in \mathbb{Z}^3 | -\frac{N}{2} \leq \bm k \leq \frac{N}{2} - 1 \Big\}.
\end{align}
The domain $\Omega$ is partitioned into a uniform grid $\Omega_N$ with mesh size $h = L/N$, defined explicitly as the set of nodes
\begin{equation}\label{set_omegaN}
\Omega_N = \left\{ (x_p, y_q, z_r) = (ph, qh, rh) \,\middle|\, (p,q,r) \in \mathcal{S}_N \right\}.
\end{equation}
All numerical variables are evaluated at the grid nodes $(x_p, y_q, z_r)$, where $x_p = ph$, $y_q = qh$ and $z_{r} = rh$ for $(p,q,r)\in \mathcal{S}_N$. Since the problem considered in this work is periodic, we define the periodic grid function space as follows:
\begin{equation}\label{set_MN}
  \mathcal{M}_N := \left\{\bm v : \Omega_N \rightarrow \mathbb{R}| v_{p+iN,q+jN,r+kN} = v_{p,q,r},\ \text{for all}\ (p,q,r)\in \mathcal{S}_N, (i,j,k)\in \mathbb{Z}^3\right\}.
\end{equation}
For any $\bm f \in \mathcal{M}_N$, let $\mu = 2\pi/L$ and $\mathrm{i} = \sqrt{-1}$. The discrete Fourier expansion is then given by
\[
f_{p,q,r} = \sum_{(k,\ell,m)\in \widehat{\mathcal{S}}_N} \hat{f}_{k,\ell,m}\mathrm{e}^{-\mathrm{i}\mu(k x_p + \ell y_q + m z_r)},\quad (p,q,r)\in \mathcal{S}_N,
\]
where the Fourier coefficients $\hat{f}_{k,\ell,m}$ are calculated via the discrete Fourier transform (DFT) as
\[
\hat{f}_{k,\ell,m} = \frac{1}{N^3} \sum_{(p,q,r)\in \mathcal{S}_N} f_{p,q,r} \mathrm{e}^{\mathrm{i}\mu(k x_p + \ell y_q + m z_r)},\quad (k,\ell,m)\in \widehat{\mathcal{S}}_N.
\]
The Fourier collocation approximations for the first- and second-order derivatives of $\bm f$ in the $x$-direction are
\[
D_x f_{p,q,r} := \sum_{(k,\ell,m)\in \widehat{S}_N}\mathrm{i}\mu k \hat{f}_{k,\ell,m} \mathrm{e}^{\mathrm{i}\mu (k x_p + \ell y_q + m z_r)},\quad D_{x}^2 f_{p,q,r} := \sum_{(k,\ell,m)\in \widehat{S}_N} -\mu^2 k^2 \hat{f}_{k, \ell, m} \mathrm{e}^{\mathrm{i}\mu (k x_p + \ell y_q + m z_r)}.
\]
The approximation for the differential operators in the $y$ and $z$ directions (i.e., $D_y$, $D_y^2$, $D_z$ and $D_z^2$) is defined analogously. Using these discrete directional derivative operators, the discrete gradient, divergence, and Laplace operators can be formulated as follows:
\[
\nabla_N \bm{f} = \left(\begin{aligned}&D_x \bm{f} \\ &D_y \bm{f} \\ &D_{z} \bm{f}\end{aligned}\right),\quad
\nabla_N\cdot\left(\begin{aligned}&\bm{f}_{1} \\ &\bm{f}_{2} \\ & \bm{f}_{3}\end{aligned}\right) = D_x \bm{f}_{1} + D_y \bm{f}_{2} + D_{z} \bm{f}_{3},\quad
\Delta_N \bm{f}  = (D_x^2 + D_y^2 + D_{z}^2) \bm{f},
\]
and it is straightforward to verify that $\nabla_N \cdot \nabla_N \bm f= \Delta_N \bm f$.

For any grid function $\bm f, \bm g \in \mathcal{M}_N$, the $\ell^2$-inner product and $\ell^2$-norm are defined, respectively, by
\[
\left<\bm f,\bm g\right> := h^3 \sum_{(p,q,r) \in \mathcal{S}_N} f_{p,q,r}\cdot g_{p,q,r} = |\Omega| \sum_{(k,\ell,m)\in \widehat{\mathcal{S}}_N}{ \hat{f}_{k,\ell,m}\cdot (\hat{g}_{k,\ell,m})^\ast}, \quad \Vert \bm g \Vert_2 := \sqrt{\left< \bm g,\bm g \right>},
\]
where {$(\cdot)^{\ast}$ denotes complex conjugation}, the second equality in the inner product definition holds by Parseval's theorem. Additionally, we introduce the $\ell^\infty$-norm and $\ell^s$-norm for $1\leq s < \infty$ as follows:
\[
\left\Vert \bm f \right\Vert _\infty := \max_{(p,q,r) \in \mathcal{S}_N} \left\vert f_{p,q,r} \right\vert,\quad \left\Vert\bm f \right\Vert _s := \Big( h^3 \sum_{(p,q,r) \in \mathcal{S}_N} |f_{p,q,r}|^s \Big)^{\frac{1}{s}}.
\]
Furthermore, given the invariance of the mass in the CH equation, it is convenient to introduce the mean-zero grid function space defined as follows:
\[
\mathring{\mathcal{M}}_{N} := \{\boldsymbol{v} \in \mathcal{M}_{N}|  \bar{\bm v} := \left< \boldsymbol{v},1 \right> = 0\} \subset \mathcal{M}_{N}.
\]
In accordance with the arguments presented in \cite{cheng2016second,cheng2019energy}, for any $\bm g \in \mathring{\mathcal{M}}_{N}$, the inverse of the discrete Laplace operator can be explicitly defined as
\[
(-\Delta_N)^{-1} g_{p,q,r} := \sum_{\substack{(k,\ell,m) \in \widehat{\mathcal{S}}_N\\ (k,\ell,m) \neq \bm 0}} (\mu^2(k^2 + \ell^2 + m^2))^{-1} \hat{g}_{k,\ell,m} \mathrm{e}^{\mathrm{i}\mu (k x_p + \ell y_q + m z_r)},\quad (p,q,r) \in \mathcal{S}_N.
\]
The following lemma presents several key inequalities for periodic grid functions, with proofs in \linkref[Appendix]{\ref{apps1}}, which will be frequently employed in the subsequent parts of this work.

\begin{lemma}\label{lem:perFun_prop}
For any $\bm f, \bm g \in \mathcal{M}_N$, the following inequalities hold
\begin{align}
  &\left\Vert \bm f \right\Vert _\infty \leq C ( |\bar{\bm f}| +  \left\Vert \nabla_N \bm f \right\Vert_2^{\frac{1}{2}} \left\Vert \Delta_N \bm f \right\Vert_2^{\frac{1}{2}}), \label{eqn:3dSob}\\[.5em]
  &\Vert \nabla_N \bm f \Vert_2 \leq C \Vert \Delta_N \bm f \Vert_2, \label{eqn:poincareIneq}\\[.5em]
  &\left\Vert \Delta_N (\bm f \bm g) \right\Vert_2 \leq C\left(\left\Vert \bm f \right\Vert _\infty \left\Vert \Delta_N \bm g \right\Vert_2 +{ \left\Vert \bm g \right\Vert_\infty \left\Vert \Delta_N \bm f \right\Vert_2}\right), \label{eqn:nonlineEst1}\\[.5em]
  &\left\Vert \nabla_N (\bm f \bm g) \right\Vert_2 \leq C\left(\left\Vert \bm f \right\Vert _\infty \left\Vert \nabla_N \bm g \right\Vert_2 + {\left\Vert \bm g \right\Vert_\infty \left\Vert \nabla_N \bm f \right\Vert_2}\right), \label{eqn:nonlineEst2}\\[.5em]
  &\left\Vert \Delta_N (\bm f^p)  \right\Vert_2 \leq C\left\Vert \bm f \right\Vert _\infty^{p-1} \left\Vert \Delta_N \bm f \right\Vert,\quad \left\Vert \nabla_N (\bm f^p)  \right\Vert_2 \leq C\left\Vert \bm f \right\Vert _\infty^{p-1} \left\Vert \nabla_N \bm f \right\Vert. \label{eqn:nonlineEst3}
\end{align}
Here and in what follows, $C$ denotes a generic positive constant that is independent of the discretization parameter but may depend on the domain $\Omega$.
\end{lemma}

\subsection{Spatial semi-discrete system}\label{Xsec4-2.2}
By introducing a constant $\kappa \geq 0$, we add and subtract the term $\kappa \Delta u$ on the right-hand side of \eqref{eqn:ch_eq}, thereby obtaining
\begin{equation}
  \label{eqn:sch_eq}
  u_t = \Delta \left( -\varepsilon^2 \Delta u + \kappa u\right) + \Delta \left( f(u) - \kappa u\right).
\end{equation}
{Applying Fourier collocation discretization in space}, we can derive the semi-discrete system for \eqref{eqn:sch_eq}, which takes the standard reaction-diffusion form:
\begin{equation}
  \label{eqn:sd_sch_eq}
  \frac{\mathrm{d}}{\mathrm{d}t}\bm u + L_\kappa \bm u = N_\kappa (\bm u),
\end{equation}
where
\begin{equation}\label{opt_LN}
  L := \varepsilon^2 \Delta_N^2,\quad  L_\kappa := L - \kappa \Delta_N, \quad N(\bm u) := \Delta_N f(\bm u) ,\quad
  N_\kappa(\bm u) := N(\bm u) - \kappa \Delta_N \bm u.
\end{equation}
The initial vector $\bm u_0$ is obtained by projecting $u_0(\bm{x})$ onto the grid $\Omega_N$. Additionally, the discrete counterpart of the energy functional \eqref{eqn:ch_energy} is given by:
\begin{equation}\label{eqn:disc_energy}
  E_N(\bm u) = \frac{\varepsilon^2}{2} \left\Vert \nabla_N \bm u \right\Vert _2^2 + \left\langle F(\bm u), 1  \right\rangle.
\end{equation}

It is straightforward to verify that the semi-discrete system \eqref{eqn:sd_sch_eq} inherits the energy dissipation law with respect to the discrete energy $E_N(\bm u)$ defined in \eqref{eqn:disc_energy}. Additionally, $E_N(\bm u)$ is inherently bounded from below, a property that guarantees the boundedness of the discrete $H^1$ semi-norm of the solutions. This result is formalized in the following lemma:

\begin{lemma}\label{lem:energy_est}
For any $\varepsilon > 0$, the discrete energy $E_N(\bm u)$ satisfies
\[
  E_N(\bm u) \geq \frac{\varepsilon^2}{2} \left\Vert \nabla_N \bm u \right\Vert _2^2 .
\]
\end{lemma}
\begin{proof}
  Recall that for $z \in \mathbb{R}$, the function $F(z) = \frac{1}{4}(z^2 - 1)^2$ is non-negative. Using this fact, we can deduce the lemma directly.
\end{proof}

\subsection{Exponential Runge--Kutta method and the energy stability results}\label{Xsec5-2.3}
For a finite $T$, let $N_t$ be a positive integer, and  $t_n := n\tau$ for $n = 0, 1, \ldots,N_t$, be the uniform partition of the time interval $[0,T]$ with $\tau = T/N_t$. The exact solution of \eqref{eqn:sd_sch_eq}  satisfies the following variation-of-constants formula
\begin{equation}
\nonumber
  \bm u(t_{n+1}) = \mathrm{e}^{-\tau L_\kappa} \bm u(t_n) + \int_{0}^{\tau} \mathrm{e}^{-(\tau - s) L_\kappa} N_\kappa (\bm u(t_{n} + s)) \mathrm{d}s.
\end{equation}
In this work, a two-stage, second-order accurate exponential Runge--Kutta scheme (ERK2) is employed for the temporal integration of the \linkref[\del{equation}\ins{Eq.}]{\eqref{eqn:sd_sch_eq}}. We denote by $\bm{u}^n$ the numerical approximation of $\bm{u}(t_n)$ and introduce an intermediate variable $\bm{u}_{n,1}$. The ERK2 scheme is then constructed by approximating $N_\kappa(\bm{u}(t_n + s))$ with $N_\kappa(\bm{u}^n)$ in the first stage and with $\tfrac{1}{2}\big(N_\kappa(\bm{u}^n) + N_\kappa(\bm{u}_{n,1})\big)$ in the second stage, i.e.,
\[
  \left\{
    \begin{aligned}
    &\bm u_{n,1} = \mathrm{e}^{-\tau L_\kappa} \bm u^n + \int_{0}^{\tau} \mathrm{e}^{-(\tau - s) L_\kappa} N_\kappa (\bm u^{n})  \mathrm{d}s,\\
    &\bm u^{n+1} = \mathrm{e}^{-\tau L_\kappa} \bm u^n + \frac{1}{2}\int_{0}^{\tau} \mathrm{e}^{-(\tau - s) L_\kappa} (N_\kappa (\bm u^{n}) + N_\kappa (\bm u_{n,1})) \mathrm{d}s.
    \end{aligned}
  \right.
\]
By evaluating the integral terms in the above expressions, we then obtain the fully discrete system based on the ERK2 scheme:
\begin{align}
 &\bm u_{n,1} = \phi_0(\tau L_\kappa) \bm u^n + \tau \phi_1(\tau L_\kappa) N_\kappa(\bm u^n)\label{eqn:etd_s1},\\
 &\bm u^{n+1} = \phi_0(\tau L_\kappa) \bm u^n + \frac{1}{2} \tau \phi_1(\tau L_\kappa) (N_\kappa(\bm u^n) + N_\kappa(\bm u_{n,1})),\label{eqn:etd_s2}
\end{align}
where the $\phi$-functions are defined by
\[
\phi_0(z) = e^{-z},\quad \phi_1(z) = \frac{1-e^{-z}}{z}.
\]
Let $\phi_i(z)$, $i = 0,1$, be functions defined on the spectrum of the linear operator $L_\kappa$ \cite{higham2008functions,li2019convergence}. Then, for any grid function $\bm{f} \in \mathcal{M}_N$ as defined in \eqref{set_MN}, $\phi_i(\tau L_\kappa) \bm f$ can be efficiently evaluated via the Fourier transform:
\begin{equation}\label{eqn:phi_operator}
  \phi_i(\tau L_\kappa) \bm f = \sum_{(k,\ell,m) \in \widehat{S}_N} \phi_i(\tau \Lambda_{k,\ell,m}) \hat{f}_{k,\ell, m} \mathrm{e}^{\mathrm{i} \mu (k x + \ell y + m z)},\quad i = 0,1.
\end{equation}
Here, $\Lambda_{k,\ell,m} =\varepsilon^2 \lambda_{k,\ell,m}^2 + \kappa \lambda_{k,\ell,m}$ and $\lambda_{k,\ell,m} = \mu^2 (k^2 + \ell^2 + m^2)$ denote the eigenvalues of $L_\kappa$ and $-\Delta_N$, respectively.

Based on the order conditions presented in \cite{hochbruck2005explicit}, the ERK2 scheme \eqref{eqn:etd_s1}--\eqref{eqn:etd_s2} possesses second-order accuracy. Consequently, the following result holds.
\begin{lemma}\label{lem:erk2_err}
Assume that the system \eqref{eqn:sd_sch_eq} has an exact solution $\bm u \in C^3([0,T], \mathcal{M}_N)$ and  $\bm u^n = \bm u(t_n)$. Then the ERK2 scheme \eqref{eqn:etd_s1}--\eqref{eqn:etd_s2} has a local truncation error of $\mathcal{O}(\tau^3)$, i.e.,
\[
  \bm u^{n+1} - \bm u(t_{n+1}) = \mathcal{O}(\tau^3).
\]
\end{lemma}
Recall that $\overline{\bm u} := \langle \bm u, 1 \rangle$ is the discrete mass of the grid function $\bm u$. The following discrete mass identities hold:
\[
  \begin{aligned}
    &\overline{\phi_0(\tau L_\kappa) \bm u^n} = \overline{\bm u^n},\ \text{ since } \phi_0(\tau \Lambda_{0,0,0}) = 1, \\
    & \overline{N_\kappa(\bm u)} = \overline{\Delta_N(f(\bm u) - \kappa \bm u)} = 0,\ \forall \bm u \in \mathcal{M}_N, \\
    & \overline{\phi_1(\tau L_\kappa)\bm v} = 0,\ \forall \bm v \in \mathring{\mathcal{M}}_{N}.
  \end{aligned}
\]
 From these identities, the following property immediately follows.
\begin{lemma}  The ERK2 scheme \eqref{eqn:etd_s1}--\eqref{eqn:etd_s2} conserves the discrete mass, that is,
  \[
    \overline{\bm u^{\,n+1}} = \overline{\bm u_{n,1}} = \overline{\bm u^n }= \overline{\bm u^0}, \quad \forall n \ge 0.
  \]
\label{Xenun9-4}
\end{lemma}
Next, we present the energy stability result of the ERK2 scheme \eqref{eqn:etd_s1}--\eqref{eqn:etd_s2} under the assumption that the numerical solution is bounded in  $\ell^{\infty}$ norm.

\begin{theorem}[Energy stability of the ERK2 scheme]\label{thm:p_stab}
Assume $\bm u^{0} \in H^{1}(\Omega)$, {and let $\{\bm u^{n},\,\bm u_{n,1},\,\bm u^{n+1}\}$ be the numerical solution computed by the ERK2 scheme \eqref{eqn:etd_s1}--\eqref{eqn:etd_s2}. Denote
\[
  M_{n} := \max\bigl\{\|\bm u^{n}\|_{\infty},\;\|\bm u_{n,1}\|_{\infty},\;\|\bm u^{n+1}\|_{\infty}\bigr\},
  \qquad
  \beta := \tfrac{3 M_{n}^{2}-1}{2}.
\]}
When {$\kappa \ge \max\left\{ 2\beta - \frac{(3-\sqrt{2}) h^{2}}{3\pi^{2}\tau}, 0 \right\}$,} it holds that
\[E_{N}(\bm u_{n,1}) \le E_{N}(\bm u^{n}),\quad E_{N}(\bm u^{n+1}) \le E_{N}(\bm u^{n}).\]
\label{Xenun1-1}
\end{theorem}

\begin{proof}
  For any $\bm u,\bm v \in \mathcal{M}_N$ satisfying $\bar{\bm u} = \bar{\bm v}$, let $M = \max(\left\Vert \bm u \right\Vert_\infty, \left\Vert \bm v \right\Vert_\infty)$ and $\beta_0 = \frac{3M^2 - 1}{2}$. For the nonlinear part in the energy functional \eqref{eqn:disc_energy}, using the Taylor expansion, we have
  \begin{equation}\label{eqn:nonlinear_est}
    \begin{aligned}
    \left\langle F(\bm v) - F(\bm u), 1 \right\rangle &
    \leq  \left\langle f(\bm u), \bm v - \bm u \right\rangle + \beta_0 \left\Vert \bm v - \bm u \right\Vert^2\\[.5em]
    &= -\left\langle (-\Delta_N)^{-1} N_\kappa (\bm u), \bm v - \bm u \right\rangle + \kappa \left\langle \bm u, \bm v - \bm u \right\rangle + \beta_0 \left\Vert \bm v - \bm u \right\Vert^2,
    \end{aligned}
  \end{equation}
 where the last equality holds by the definition of the nonlinear operator $N_\kappa(\cdot)$ in  \eqref{opt_LN}.
In addition, by using the algebraic identity $a^2 - b^2 = 2a(a-b) - (a-b)^2$, the first term in the discrete energy \eqref{eqn:disc_energy} satisfies
  \begin{equation}\label{eqn:linear_est}
    \begin{aligned}
    \frac{\varepsilon^2}{2}\left\langle \nabla_N \bm v, \nabla_N \bm v \right\rangle - \frac{\varepsilon^2}{2} \left\langle \nabla_N \bm u, \nabla_N \bm u \right\rangle =& \varepsilon^2 \left\langle \nabla_N \bm v, \nabla_N (\bm v - \bm u) \right\rangle - \frac{\varepsilon^2}{2}\left\Vert \nabla_N (\bm v - \bm u) \right\Vert_2^2\\
    = & \left\langle  (- \Delta_N)^{-1} L_\kappa \bm v, \bm v - \bm u \right\rangle - \kappa \left\langle \bm v, \bm v - \bm u \right\rangle - \frac{\varepsilon^2}{2}\left\Vert \nabla_N (\bm v - \bm u) \right\Vert_2^2.
    \end{aligned}
  \end{equation}
Therefore, by combining \eqref{eqn:nonlinear_est} and \eqref{eqn:linear_est}, we obtain
  \begin{equation}\label{eqn:energ_diff}
  \begin{aligned}
   & E_N(\bm v) - E_N(\bm u) \leq  -\left\langle  (- \Delta_N)^{-1} (-L_\kappa \bm v + N_\kappa (\bm u)), \bm v - \bm u \right\rangle - (\kappa - \beta_0) \left\Vert \bm v- \bm u \right\Vert_2^2 - \frac{\varepsilon^2}{2}\left\Vert \nabla_N (\bm v - \bm u) \right\Vert_2^2\\
    \leq & -\left\langle  (- \Delta_N)^{-1} (-L_\kappa \bm v + N_\kappa (\bm u)), \bm v - \bm u \right\rangle- \frac{1}{2} \left\langle  (- \Delta_N)^{-1} L_\kappa( \bm v -  \bm u ), \bm v - \bm u \right\rangle - (\frac{\kappa}{2} - \beta_0) \left\Vert \bm v- \bm u \right\Vert_2^2.
    \end{aligned}
  \end{equation}
Rearranging the formulation for the intermediate stage in \eqref{eqn:etd_s1}, we have
  \begin{equation}\label{eqn:1st_stage_re}
    \begin{aligned}
     &-L_\kappa \bm u_{n,1} + N_\kappa(\bm u^n) \\
    =& -L_\kappa \bm u_{n,1} + (\tau \phi_1(\tau L_\kappa))^{-1} (\bm u_{n,1} - \phi_0(\tau L_\kappa)\bm u^{n})\\
    =& -L_\kappa (\bm u_{n,1} - \bm u^{n}) + (\tau \phi_1(\tau L_\kappa))^{-1} (\bm u_{n,1} - \bm u^{n})\\
    =& \frac{1}{\tau} S_1(\tau L_\kappa) (\bm u_{n,1} - \bm u^{n}),
    \end{aligned}
  \end{equation}
where $S_1(z) = \frac{z}{e^z - 1}$. Substituting \eqref{eqn:1st_stage_re} into \eqref{eqn:energ_diff} yields
  \[
  \begin{aligned}
    E_N(\bm u_{n,1}) - E_N(\bm u^n) \leq & -\frac{1}{\tau}\left\langle (-\Delta_N)^{-1} \left(S_1(\tau L_\kappa) + \frac{\tau}{2} L_\kappa\right) (\bm u_{n,1} - \bm u^n), \bm u_{n,1} - \bm u^n \right\rangle - (\frac{\kappa}{2} - \beta) \left\Vert \bm u_{n,1} - \bm u^n \right\Vert_2^2\\
    \leq & - (\frac{1}{\Lambda_N\tau} + \frac{\kappa}{2} - \beta)  \Vert \bm u_{n,1} - \bm u^n \Vert^2,
  \end{aligned}
  \]
  where $\Lambda_{N} = \frac{3\pi^2}{h^2}$ denotes the largest eigenvalue of $-\Delta_N$. This result indicates that if the stabilization parameter satisfies
{  \[
    \kappa \geq \max\left\{  2 \beta - \frac{2h^2}{3\pi^2 \tau}, 0\right\},
  \]}
  the discrete energy corresponding to the first stage solution satisfies the strict monotonic decay property $E_N(\boldsymbol{u}_{n,1}) \leq E_N(\boldsymbol{u}^n)$.
Similarly, from \eqref{eqn:etd_s2}, we can obtain
\begin{equation}\label{eqn:2nd_stage_re}
    \begin{aligned}
     &-L_\kappa \bm u^{n+1} + N_\kappa(\bm u_{n,1})\\
    =& -L_\kappa \bm u^{n+1} + \left(\frac{1}{2}\tau \phi_1(\tau L_\kappa)\right)^{-1}\left(\bm u^{n+1} - \phi_0(\tau L_\kappa) \bm u^n - \frac{1}{2}\tau \phi_1(\tau L_\kappa)N_\kappa (\bm u^n)\right)\\
    =& -L_\kappa \bm u^{n+1} + \left(\frac{1}{2}\tau \phi_1(\tau L_\kappa)\right)^{-1}\left(\bm u^{n+1} - \phi_0(\tau L_\kappa) \bm u^n - \frac{1}{2} \left(\bm u_{n,1} - \phi_0(\tau L_\kappa)\bm u^{n}\right)\right)\\
    =& -L_\kappa \bm u^{n+1} + \left(\frac{1}{2}\tau \phi_1(\tau L_\kappa)\right)^{-1}\left(\bm u^{n+1} - \frac{1}{2}\bm u_{n,1} - \frac{1}{2}\phi_0(\tau L_\kappa) \bm u^n\right)\\
    =& -L_\kappa (\bm u^{n+1} - \bm u^{n}) + \left(\frac{1}{2}\tau \phi_1(\tau L_\kappa)\right)^{-1} (\bm u^{n+1} - \bm u_{n,1}) + \left(\tau \phi_1(\tau L_\kappa)\right)^{-1}(\bm u_{n,1} - \bm u^n) \\
    =& \frac{1}{\tau} S_2(\tau L_\kappa) (\bm u^{n+1} - \bm u_{n,1}) + \frac{1}{\tau} S_1(\tau L_\kappa) (\bm u_{n,1} - \bm u^{n}),
    \end{aligned}
\end{equation}
where $S_2(z) = \frac{z(e^z+1)}{e^z-1}$. It then follows that
\[
    \begin{aligned}
       E_N(\bm u^{n+1}) - E_N(\bm u^n) =& E_N(\bm u^{n+1}) - E_N(\bm u_{n,1}) + E_N(\bm u_{n,1}) - E_N(\bm u^n) \\
      \leq & -\frac{1}{\tau} \left\langle  (-\Delta_N)^{-1} \left(S_2(\tau L_\kappa) + \frac{\tau}{2} L_\kappa\right) (\bm u^{n+1} - \bm u_{n,1}), \bm u^{n+1} - \bm u_{n,1} \right\rangle\\
      &-\frac{1}{\tau} \left\langle (-\Delta_N)^{-1} S_1(\tau L_\kappa) (\bm u_{n,1} - \bm u^{n}), \bm u^{n+1} - \bm u_{n,1} \right\rangle  \\
      &- \frac{1}{\tau} \left\langle (-\Delta_N)^{-1}\left(S_1(\tau L_\kappa) +\frac{\tau}{2} L_\kappa\right)  (\bm u_{n,1} - \bm u^n), \bm u_{n,1} - \bm u^n \right\rangle \\
      &- \left(\frac{\kappa}{2} - \beta\right) \left\Vert \bm u_{n,1} - \bm u^n \right\Vert_2^2 - \left(\frac{\kappa}{2} - \beta\right) \left\Vert \bm u^{n+1} - \bm u_{n,1} \right\Vert_2^2.
    \end{aligned}
\]
{To express the energy difference as a quadratic form of a symmetric matrix, we first define the increment vectors
\[
\delta U = \left(\begin{array}{c}
  \bm u_{n,1} - \bm  u^{n}\\
  \bm u^{n+1} - \bm u_{n,1}
\end{array}\right),\quad
(-\bm{ \Delta_N})^{-1} = \left[\begin{array}{cc}
  (- \Delta_N)^{-1} & 0\\
  0 & (- \Delta_N)^{-1}
\end{array}\right],
\]
and the symmetric block operator matrix
\[
    \Phi(\tau L_\kappa) = \left[
    \begin{array}{cc}
      S_1(\tau L_\kappa) + \frac{1}{2}\tau L_\kappa &  \frac{1}{2}S_1(\tau L_\kappa)  \\
      \frac{1}{2}S_1(\tau L_\kappa)  & S_2(\tau L_\kappa) + \frac{1}{2}\tau L_\kappa
    \end{array}
    \right].
\]
Then, the energy difference can be rewritten as
\[
\begin{aligned}
E_N(\boldsymbol{u}^{n+1}) - E_N(\boldsymbol{u}^n)
\leq &- \frac{1}{\tau}\left\langle \delta U, (-\bm \Delta_N)^{-1} \Phi(\tau L_\kappa) \delta U\right\rangle - \left(\frac{\kappa}{2} - \beta\right) \langle \delta U, \delta U \rangle.
\end{aligned}
\]
For any $z\geq0$, the eigenvalues of the scalar matrix $\Phi(z)$ are
\[
\lambda(\Phi(z)) = \frac{1}{2}\left[3S_1(z)+ 2z \pm \sqrt{2S_1(z)^2 +2zS_1(z)+z^2}\right].
\]
The minimum eigenvalue is strictly positive and uniformly bounded below by $\frac{3 - \sqrt{2}}{2} > 0$ for all $z\geq0$, thus, we have
\[
\left\langle \delta \boldsymbol{U}, (-\bm{\Delta})^{-1}\Phi(\tau L_\kappa) \delta \boldsymbol{U} \right\rangle \geq \frac{3-\sqrt{2}}{2\Lambda_N} \left\langle \delta \boldsymbol{U}, \delta \boldsymbol{U} \right\rangle.
\]
Substituting into the energy difference inequality yields
\[
E_N(\boldsymbol{u}^{n+1}) - E_N(\boldsymbol{u}^n) \leq  - \left(\frac{3-\sqrt{2}}{2\Lambda_N \tau} + \frac{\kappa}{2} - \beta\right) \left\langle  \delta \boldsymbol{U}, \delta \boldsymbol{U} \right\rangle.
\]
For the 3D uniform grid with step size $h$, substituting $\Lambda_N = \frac{3\pi^2}{h^2}$ gives
\[
E_N(\boldsymbol{u}^{n+1}) - E_N(\boldsymbol{u}^n) \leq  - \left(\frac{(3-\sqrt{2})h^2}{6\pi^2 \tau} + \frac{\kappa}{2} - \beta \right) \left\langle  \delta \boldsymbol{U}, \delta \boldsymbol{U} \right\rangle.
\]
Thus, we conclude that when the stabilization parameter satisfies
\[
\kappa \geq \max\left\{  2 \beta - \frac{(3-\sqrt{2})h^2}{3\pi^2 \tau}, 0\right\},
\]
the discrete energy satisfies the strict monotonic decay property $E_N(\boldsymbol{u}^{n+1}) \leq E_N(\boldsymbol{u}^n)$. Consequently, this completes the proof and yields the desired results.
}
\end{proof}

{In phase-field simulations with a double-well potential, the order parameter is typically observed to remain of $\mathcal{O}(1)$ in the maximum norm. This observation makes \linkref[Theorem]{\ref{thm:p_stab}} useful as a practical energy-dissipation criterion: once a maximum-norm bound for the numerical solution and the stage values is available, the stabilization parameter $\kappa$ can be chosen or verified accordingly. From the analytical viewpoint, however, such a bound should not be imposed as an a priori assumption. Although the continuous CH equation admits uniform-in-time $L^\infty$ bounds under suitable regularity assumptions \cite{elliott1986cahn,novick2008cahn,temam2012infinite}, the corresponding fully discrete maximum-norm bound requires a separate proof. The theorem below provides this missing ingredient: when the time step is sufficiently small, with a threshold independent of the spatial discretization parameter, the numerical solution admits uniform-in-time $H^1$, $H^2$, and consequently $\ell^\infty$ bounds.}

\begin{theorem}\label{thm:long_time_stability}
Suppose the initial data for the scheme \eqref{eqn:etd_s1}--\eqref{eqn:etd_s2} satisfy $u_0(\bm{x}) \in H^4(\Omega)$, and let the corresponding initial energy satisfy $E(u_0(\bm{x})) \le C_e$. Then, for any $\kappa \ge 0$ and under a time-step constraint $\tau \le \tau_s{:= \min\{\tau_1, \tau_2\}}$, {with $\tau_1$ and $\tau_2$ given in \eqref{eqn:stepConst_1} and \eqref{eqn:stepConst_2}, respectively,} there exist constants $\widehat{C}_{\varepsilon,0}$ and $\widehat{C}_{\varepsilon,1}$, independent of the final time $T$, and a constant $\ddot{C}_0$ depending on the initial data and the regularity assumption. For all $n \ge 0$, the numerical solution satisfies the following bounds:
    \[
    \begin{aligned}
        \Vert \nabla_N \bm u^{n}\Vert_2^2 \leq 2\varepsilon^{-2} {C}_e,
        \quad \Vert \Delta_N \bm u^{n}\Vert_2^2 \leq \Big(1 - \frac{\varepsilon^2 \tau}{\widehat{C}_{\varepsilon,0} +  \varepsilon^2 \tau}\Big)^{n} \ddot{C}_0^2 + C \varepsilon^{-2}\widehat{C}_{\varepsilon,1},
    \end{aligned}
    \]
    and the $\ell^\infty$ bound for the numerical solution is consequently obtained via the Sobolev embedding inequality \eqref{eqn:3dSob}.
\end{theorem}

{\opecolor{blue}
\begin{remark}  While the same stabilization parameter $\kappa$ is used in \linkref[Theorems]{\ref{thm:p_stab}} and~\ref{thm:long_time_stability}, they serve distinct purposes. \linkref[Theorem]{\ref{thm:p_stab}} provides a one-step energy dissipation criterion for the ERK2 scheme: assuming that the relevant stage values admit an $\ell^\infty$ bound, energy dissipation follows when $\kappa$ satisfies the stated lower bound depending on this maximum norm and on $\tau$. \linkref[Theorem]{\ref{thm:long_time_stability}} proves that for any fixed $\kappa \ge 0$, the numerical solution admits uniform-in-time $H^1$, $H^2$, and consequently $\ell^\infty$ bounds under the time-step condition $\tau \le \tau_s := \min\{\tau_1, \tau_2\}$. This condition is needed to close the long-time bootstrap argument for the maximum-norm bound. The threshold $\tau_s$ is independent of the spatial mesh size $N$ but may depend on $\varepsilon$, $\kappa$, and constants related to the initial data and regularity assumptions. Once this uniform maximum-norm bound is available, it can be used in \linkref[Theorem]{\ref{thm:p_stab}} to formulate the corresponding energy-dissipation condition on the same parameter $\kappa$.
\label{Xenun4-1}
\end{remark}
}

\section{The proof of \linkref[Theorem]{\ref{thm:long_time_stability}}: long-time numerical stability analysis}\label{Xsec6-3}\label{sec:stab}
In this section, we focus on deriving uniform-in-time $H^m$ estimates ($m = 1, 2$) for the numerical solution at all stages, and then obtain $\ell^\infty$ boundedness via the discrete Sobolev embedding inequality \eqref{eqn:3dSob}. Motivated by the long-time stability analysis techniques proposed in \cite{gottlieb2012long,cheng2016long}, our goal is to derive estimates independent of the final time $T$. To this end, we first carry out preparatory work, including setting up a suitable a priori assumption and establishing preliminary estimates.
\subsection{Some preliminary estimates and setup}\label{Xsec7-3.1}
For a better analysis of the interaction between the linear and nonlinear terms, we first recast the ERK2 scheme \eqref{eqn:etd_s1}--\eqref{eqn:etd_s2} into the following equivalent form:
\begin{align}
  &\bm u_{n,1} - \bm u^{n} + \tau \phi_1(\tau L_\kappa) L_\kappa \bm u^n = \tau \phi_1(\tau L_\kappa) N_\kappa (\bm u^n), \label{eqn:rfd_1}\\
  &\bm u^{n+1} - \bm u^{n} + \tau \phi_1(\tau L_\kappa) L_\kappa \bm u^n = \frac{1}{2}\tau \phi_1(\tau L_\kappa)(N_\kappa(\bm u^n) + N_\kappa(\bm u_{n,1})). \label{eqn:rfd_2}
\end{align}

Furthermore, the following algebraic bounds for the $\phi$-functions are introduced, which are summarized in the lemma below.
\begin{lemma}For the function $\phi_1(z) = \frac{1-e^{-z}}{z}$, the following inequalities are valid:
\begin{equation}\label{eqn:phi_esitmate}
   0 \leq \phi_1(z) \leq 1,\quad 1 \leq (\phi_1(z))^{-1} \leq 1+z,\quad \forall z \geq 0,
\end{equation}
and
\begin{equation}
  (\phi_1(z))^{-1} \geq z,\quad \forall z \in \mathbb{R}.
\label{Xeqn18-3.4}
\end{equation}
\label{Xenun10-5}
\end{lemma}

To avoid excessive complexity of the formulas during the stability analysis, we introduce the following linear operators:
\begin{equation}\label{eqn:linear_operators}
  G_N^{(1)} = (-\phi_1(\tau L_\kappa)L_\kappa \Delta_N)^{\frac{1}{2}},\quad
  G_N^{(2)} = (\phi_1(\tau L_\kappa)L_\kappa \Delta_N^2)^{\frac{1}{2}},\quad
  \tilde{G}_N^{(2)} = (L_\kappa \Delta_N^2)^{\frac{1}{2}}.
\end{equation}
These operators are defined in the same fashion as \eqref{eqn:phi_operator}. In more detail, for any grid function $\bm f\in \mathcal{M}_N$, we have
\begin{align*}
& G_N^{(1)}\bm f  = \sum_{(k,\ell,m) \in \widehat{\mathcal{S}}_N} (\phi_1(\tau \Lambda_{k,\ell,m})\Lambda_{k,\ell,m} \lambda_{k,\ell,m} )^{\frac{1}{2}}\hat{f}_{k,\ell, m} \mathrm{e}^{\mathrm{i} \mu (k x + \ell y + m z)},\\
& G_N^{(2)}\bm f  = \sum_{(k,\ell,m) \in \widehat{\mathcal{S}}_N} (\phi_1(\tau \Lambda_{k,\ell,m})\Lambda_{k,\ell,m} \lambda_{k,\ell,m}^2 )^{\frac{1}{2}}\hat{f}_{k,\ell, m} \mathrm{e}^{\mathrm{i} \mu (k x + \ell y + m z)},\\
& \tilde{G}_N^{(2)}\bm f  = \sum_{(k,\ell,m) \in \widehat{\mathcal{S}}_N} (\Lambda_{k,\ell,m} \lambda_{k,\ell,m}^2 )^{\frac{1}{2}}\hat{f}_{k,\ell, m} \mathrm{e}^{\mathrm{i} \mu (k x + \ell y + m z)}.
\end{align*}
It follows directly from these definitions that the operators are positive semi-definite and satisfy
\[
  \left\langle \phi_1(\tau L_\kappa) L_\kappa \bm f, -\Delta_N \bm f \right\rangle = \Vert G_N^{(1)} \bm f \Vert_2^2,\quad \left\langle \phi_1(\tau L_\kappa) L_\kappa \bm f, \Delta_N^2 \bm f \right\rangle = \Vert G_N^{(2)} \bm f \Vert_2^2,\quad \left\langle L_\kappa \bm f, \Delta_N^2 \bm f \right\rangle = \Vert \tilde{G}_N^{(2)} \bm f \Vert_2^2.
\]
Moreover, to deal with the nonlinear structure of the CH \linkref[\del{equation}\ins{Eq.}]{\eqref{eqn:ch_eq}}, we introduce the inverse operator
\[
  (\phi_1(\tau L_\kappa))^{-1} \bm f = (1-e^{-\tau L_\kappa})^{-1}(\tau L_\kappa) \bm f = \sum_{(k,\ell,m) \in \widehat{\mathcal{S}}_N} \frac{\tau \Lambda_{k,\ell,m}}{1 - e^{-\tau \Lambda_{k,\ell,m}}}\hat{f}_{k,\ell, m} \mathrm{e}^{\mathrm{i} \mu (k x + \ell y + m z)},
\]
{\opecolor{blue}
and
\[
  (\phi_1(\tau L_\kappa))^{-\frac12} \bm f = \sum_{(k,\ell,m) \in \widehat{\mathcal{S}}_N} \left(\frac{\tau \Lambda_{k,\ell,m}}{1 - e^{-\tau \Lambda_{k,\ell,m}}}\right)^{\frac12} \hat{f}_{k,\ell, m} \mathrm{e}^{\mathrm{i} \mu (k x + \ell y + m z)}.
\]}
Accordingly, the corresponding summation-by-parts identity reads
\[
  \left\langle (\phi_1(\tau L_\kappa))^{-1} \bm f, \bm f \right\rangle = \Vert (\phi_1(\tau L_\kappa))^{-\frac{1}{2}} \bm f \Vert_2^2.
\]
The following preliminary estimates will be required in the subsequent analysis.
\begin{lemma}\label{lem:operator_est}
For any grid functions $\bm f \in \mathcal{M}_N$, we have
\begin{align}\label{eqn:diffusion_estimate}
&\Vert \nabla_N \bm f \Vert_2^2 \geq \tau \Vert G_N^{(1)}\bm f \Vert_2^2,\quad \Vert \Delta_N \bm f \Vert_2^2 \geq \tau \Vert G_N^{(2)} \bm f \Vert_2^2,\quad \Vert\left( \phi_1(\tau L_\kappa)\right)^{-\frac{1}{2}}\Delta_N \bm f \Vert_2^2 \geq \tau \Vert \tilde{G}_N^{(2)}\bm f \Vert_2^2,\\
  \label{eqn:diffusion_split_1}
&  \Vert G_N^{(1)}\bm f \Vert_2^2 = \varepsilon^2 \Vert \left( \phi_1(\tau L_\kappa)\right)^{\frac{1}{2}} \Delta_N \nabla_N \bm f \Vert_2^2 + \kappa \Vert \left( \phi_1(\tau L_\kappa)\right)^{\frac{1}{2}} \Delta_N \bm f \Vert_2^2,\\
  \label{eqn:diffusion_split_2}
&  \Vert G_N^{(2)} \bm f \Vert_2^2 = \varepsilon^2 \Vert \left( \phi_1(\tau L_\kappa)\right)^{\frac{1}{2}} \Delta_N^2 \bm f \Vert_2^2 + \kappa \Vert \left( \phi_1(\tau L_\kappa)\right)^{\frac{1}{2}} \Delta_N \nabla_N \bm f \Vert_2^2, \\
  \label{eqn:diffusion_split_3}
&  \Vert \tilde{G}_N^{(2)}\bm f \Vert_2^2 = \varepsilon^2 \Vert \Delta_N^2 \bm f \Vert_2^2 + \kappa \Vert \Delta_N \nabla_N \bm f \Vert_2^2,\\
  \label{eqn:diffusion_split_4}
&   \Vert \Delta_N \bm f \Vert_2^2 \leq  \Vert \left(\phi_1(\tau L_\kappa)\right)^{-\frac{1}{2}} \Delta_N \bm f \Vert_2^2 \leq \Vert \Delta_N \bm f \Vert_2^2 + \tau (\varepsilon^2 \Vert \Delta_N^2 \bm f\Vert_2^2 + \kappa \Vert \Delta_N \nabla_N \bm f \Vert_2^2 ).
\end{align}
\end{lemma}
\begin{proof}
For any $\bm f \in \mathcal{M}_N$, an application of Parseval's equality yields
\[
    \Vert G_{N}^{(1)} \bm f \Vert_2^2 = |\Omega|\sum_{(k,\ell,m)\in \widehat{\mathcal{S}}_N} \frac{1 - \mathrm{e}^{-\tau \Lambda_{k,\ell,m}}}{\tau \Lambda_{k,\ell,m}} \Lambda_{k,\ell,m} \lambda_{k,\ell,m} |\hat{f}_{k,\ell,m}|^2.
\]
By noting the trivial fact that $1 - \mathrm{e}^{-\tau \Lambda_{k,\ell,m}} \leq 1$ for $\tau \Lambda_{k,\ell,m} \geq 0$,  we have
\[
    \tau \Vert G_{N}^{(1)} \bm f \Vert_2^2 = |\Omega|\sum_{(k,\ell,m)\in \widehat{\mathcal{S}}_N} (1 - \mathrm{e}^{-\tau \Lambda_{k,\ell,m}}) \lambda_{k,\ell,m} |\hat{f}_{k,\ell,m}|^2 \leq |\Omega|\sum_{(k,\ell,m)\in \widehat{\mathcal{S}}_N} \lambda_{k,\ell,m} |\hat{f}_{k,\ell,m}|^2 = \Vert \nabla_N \bm f \Vert_2^2.
\]
Thus, the proof of the first inequality in \eqref{eqn:diffusion_estimate} is completed. The second and third inequalities can be derived by following the same line of reasoning.

In addition, given the definition of the linear operator $L_\kappa$ in \eqref{opt_LN}, identity \eqref{eqn:diffusion_split_1} can be directly derived by applying the summation by parts formula:
\[
\begin{aligned}
    \Vert G_{N}^{(1)} \bm f \Vert_2^2 = \langle \phi_1(\tau L_\kappa) L_\kappa \bm f, -\Delta_N \bm f \rangle &= -\varepsilon^2 \langle \phi_1(\tau L_\kappa) \Delta_N^2 \bm f, \Delta_N \bm f \rangle + \kappa \langle \phi_1(\tau L_\kappa) \Delta_N \bm f, \Delta_N \bm f \rangle \\
    &= \varepsilon^2 \Vert \left( \phi_1(\tau L_\kappa)\right)^{\frac{1}{2}} \Delta_N \nabla_N \bm f \Vert_2^2 + \kappa \Vert \left( \phi_1(\tau L_\kappa)\right)^{\frac{1}{2}} \Delta_N \bm f \Vert_2^2.
\end{aligned}
\]
Identities \eqref{eqn:diffusion_split_2} and \eqref{eqn:diffusion_split_3} can be derived in the same manner. The estimate in \eqref{eqn:diffusion_split_4} can be obtained directly by using the result in \eqref{eqn:phi_esitmate}. The proof is completed.
\end{proof}

In the remainder of this section, we first establish a long-time stability analysis of the numerical solution in the $H^1$ and $H^2$ norms via the energy technique. By applying the Sobolev inequality \eqref{eqn:3dSob}, we obtain a uniform-in-time bound for the numerical solution in  $\ell^\infty$ norm. This result, in turn, leads to a uniform-in-time stability estimate for the ERK2 scheme \eqref{eqn:etd_s1}--\eqref{eqn:etd_s2} and completes the proof of \linkref[Theorem]{\ref{thm:long_time_stability}}.

To carry out uniform-in-time estimates for the numerical solution, we first introduce a key a priori assumption regarding energy dissipation at the previous time step:
\begin{equation}\label{eqn:apriori_En}
E(\bm u^n) \leq E(\bm u^0) := C_e,
\end{equation}
which will be recovered at the next time step. By leveraging this discrete energy assumption together with \linkref[Lemma]{\ref{lem:energy_est}}, the following $H^1$ norm bound for the numerical solution can be derived:
\begin{equation}\label{eqn:apriori_H1}
\left\Vert \nabla_N \bm u^n  \right\Vert_2 \leq \dot{C}_{n} :=  \varepsilon^{-1} (2 C_e)^{\frac{1}{2}}.
\end{equation}
Additionally, since the ERK2 scheme preserves mass at the discrete level, we impose the zero average condition on the initial data $\bm u^{0}$ to simplify the subsequent analysis. This leads to
\begin{equation}
    \overline{\bm u^{n+1}} = \overline{\bm u_{n,1}} = \overline{\bm u^{n}}= 0.
\label{Xeqn27-3.13}
\end{equation}

On the other hand, to facilitate the uniform-in-time $H^2$ estimate of the numerical solution, a regularity assumption is imposed on the initial data:
\[
\Vert (\phi{_1}(\tau L_\kappa))^{-\frac{1}{2}}\Delta_N \bm u^{0}  \Vert_2 \leq \ddot{C}_0.
\]
According to the preliminary estimate \eqref{eqn:diffusion_split_4}, this assumption essentially imposes an $H^2$ regularity requirement on the initial data. Notably, $H^4$ regularity is sufficient to theoretically guarantee the boundedness of this quantity. To facilitate the derivation of the $H^2$ estimate for the numerical solution, we impose the following a priori assumption at the preceding time level
\begin{equation}\label{eqn:apriori_psi}
\left\Vert (\phi_1(\tau L_\kappa))^{-\frac{1}{2}}\Delta_N \bm u^{n} \right\Vert_2^2 \leq \ddot{C}_n^2 := \Big(1 - \frac{\varepsilon^2 \tau}{\widehat{C}_{\varepsilon,0} +  \varepsilon^2 \tau}\Big)^{n} \ddot{C}_0^2 + C \Big(1- \Big(1 - \frac{\varepsilon^2 \tau}{\widehat{C}_{\varepsilon,0} + \varepsilon^2 \tau}\Big)^{n}\Big) (1 + \varepsilon^2\tau ) \varepsilon^{-2}\widehat{C}_{\varepsilon,1},
\end{equation}
where $\widehat{C}_{\varepsilon,0}, \widehat{C}_{\varepsilon,1}$ are global-in-time constants, independent of the final time $T$, which may depend on  $\varepsilon^{-1}$ and $\kappa$ in a polynomial pattern, and their specific values will be determined in subsequent analysis.

The first term in the expansion of $\ddot{C}_n$ in \eqref{eqn:apriori_psi} represents an exponential decay of the contribution from the initial data, reflecting the parabolic nature of the CH \linkref[\del{equation}\ins{Eq.}]{\eqref{eqn:ch_eq}}, while the second term has a uniform-in-time bound under the assumption $\tau \varepsilon^2 \leq 1$. Consequently, an  $H^2$ bound for the numerical solution at $t_n$ becomes available
\begin{equation}\label{eqn:apriori_H2}
\left\Vert \Delta_N \bm u^n \right\Vert_2 \leq \left\Vert (\phi_1(\tau L_\kappa))^{-\frac{1}{2}}\Delta_N \bm u^{n} \right\Vert_2 \leq \ddot{C}_{n},
\end{equation}
in which the first inequality comes from the inequality \eqref{eqn:phi_esitmate}.  Moreover, by combining the a priori estimates in \eqref{eqn:apriori_H1} and \eqref{eqn:apriori_H2} and applying the Sobolev inequality \eqref{eqn:3dSob}, we obtain an explicit upper bound for the maximum norm
\[
\left\Vert \bm u^n \right\Vert_\infty  \leq \tilde{C}_n = C \dot{C}_{n}^{\frac{1}{2}}\ddot{C}_{n}^{\frac{1}{2}},
\]
which will be used in the subsequent nonlinear analysis.

\subsection{Preliminary estimate for \texorpdfstring{$\Vert\nabla_N \bm u_{n,1}\Vert_2$}{}}\label{Xsec8-3.2}
Taking the discrete $\ell^2$ inner product of \eqref{eqn:rfd_1} with $-2 \Delta_N \bm u_{n,1}$ yields
\begin{equation}\label{eqn:rough1_h1_1}
  \left\langle \bm u_{n,1} - \bm u^n , -2\Delta_N \bm u_{n,1} \right\rangle - 2\tau \left\langle \phi_1(\tau L_\kappa) L_\kappa \bm u^n, \Delta_N \bm u_{n,1} \right\rangle = -2 \tau \left\langle \phi_1(\tau L_\kappa) N_\kappa(\bm u^n), \Delta_N \bm u_{n,1} \right\rangle.
\end{equation}
For the first term on the left-hand side, an application of the summation-by-parts formula gives
\begin{equation}\label{eqn:rough1_h1_2}
  \begin{aligned}
  \left\langle \bm u_{n,1} - \bm u^n , -2\Delta_N \bm u_{n,1} \right\rangle =& \left\langle \nabla_N \bm u_{n,1} - \nabla_N \bm u^n , 2\nabla_N \bm u_{n,1} \right\rangle \\
   =& \left\Vert \nabla_N \bm u_{n,1} \right\Vert_2^2 - \left\Vert \nabla_N \bm u^n \right\Vert_2^2 + \left\Vert \nabla_N (\bm u_{n,1} - \bm u^n) \right\Vert_2^2.
  \end{aligned}
\end{equation}
Next, consider the linear diffusion term in \eqref{eqn:rough1_h1_1}. Through direct computation, we obtain
\begin{equation}\label{eqn:rough1_h1_3}
  \begin{aligned}
  -2 \left\langle \phi_1(\tau L_\kappa) L_\kappa \bm u^n, \Delta_N \bm u_{n,1} \right\rangle =& 2\left\langle G_N^{(1)}\bm u^n, G_N^{(1)}\bm u_{n,1} \right\rangle \\
  =& \Vert G_N^{(1)}\bm u^n \Vert_2^2 + \Vert G_N^{(1)}\bm u_{n,1} \Vert_2^2 - \Vert G_N^{(1)} (\bm u^n - \bm u_{n,1}) \Vert_2^2.
  \end{aligned}
\end{equation}
Moreover, it follows from \eqref{eqn:diffusion_estimate} that
\begin{equation}\label{eqn:rough1_h1_4}
\lVert \nabla_{N} (\bm u^n - \bm u_{n,1}) \rVert_{2}^2 - \tau \lVert G_N^{(1)} (\bm u^n - \bm u_{n,1}) \rVert_{2}^2 \geq 0.
\end{equation}
Substituting \eqref{eqn:rough1_h1_2}-\eqref{eqn:rough1_h1_4} into \eqref{eqn:rough1_h1_1}, we obtain
\begin{equation}\label{eqn:rough1_h1_5}
  \left\Vert \nabla_N \bm u_{n,1} \right\Vert_2^2 - \left\Vert \nabla_N \bm u^n \right\Vert_2^2 + \tau (\Vert G_N^{(1)} \bm u^n \Vert_2^2 + \Vert G_N^{(1)} \bm u_{n,1} \Vert_2^2 ) \leq -2 \tau \left\langle \phi_1(\tau L_\kappa) N_\kappa(\bm u^n), \Delta_N \bm u_{n,1} \right\rangle.
\end{equation}
Here, using the estimate of operator $G_N^{(1)}$ provided in \eqref{eqn:diffusion_split_1}, we have
\begin{equation}
\label{ineq_G1}
\begin{aligned}
    \Vert G_N^{(1)} \bm u^n \Vert_2^2 + \Vert G_N^{(1)} \bm u_{n,1} \Vert_2^2
    =& \varepsilon^2 (\Vert \left( \phi_1(\tau L_\kappa)\right)^{\frac{1}{2}} \Delta_N \nabla_N \bm u^n \Vert_2^2 + \Vert \left( \phi_1(\tau L_\kappa)\right)^{\frac{1}{2}} \Delta_N \nabla_N \bm u_{n,1} \Vert_2^2) \\ &+ \kappa (\Vert \left( \phi_1(\tau L_\kappa)\right)^{\frac{1}{2}} \Delta_N \bm u^n \Vert_2^2  + \Vert \left( \phi_1(\tau L_\kappa)\right)^{\frac{1}{2}} \Delta_N \bm u_{n,1} \Vert_2^2).
\end{aligned}
\end{equation}
By applying \eqref{opt_LN}, the nonlinear term on the right-hand side of \eqref{eqn:rough1_h1_5} becomes
\begin{equation}\label{nonl_inn}
  \begin{aligned}
&- 2\left\langle \phi_1(\tau L_\kappa) N_\kappa(\bm u^n), \Delta_N \bm u_{n,1} \right\rangle\\
=& -2\left< \phi_{1}(\tau L_{\kappa})\Delta_{N} (\bm u^n)^3, \Delta_{N}\bm u_{n,1} \right> + 2(\kappa + 1) \left< \phi_{1}(\tau L_{\kappa}) \Delta_{N} \bm u^n , \Delta_{N}\bm u_{n,1} \right>.
  \end{aligned}
\end{equation}
Using \eqref{eqn:nonlineEst3} and \eqref{eqn:phi_esitmate}, the first nonlinear inner product can be bounded as
\begin{equation}\label{eqn:rough1_h1_7}
\begin{aligned}
-2 \left< \phi_{1}(\tau L_{\kappa})\Delta_{N} (\bm u^n)^3, \Delta_{N}\bm u_{n,1} \right> \leq & \frac{2}{\varepsilon^2}\lVert \left(\phi_{1}(\tau L_{\kappa})\right)^{\frac{1}{2}} \nabla_{N} (\bm u^n)^3 \rVert_{2}^2 + \frac{\varepsilon^2}{2}\lVert \left(\phi_{1}(\tau L_{\kappa})\right)^{\frac{1}{2}}\nabla_{N} \Delta_{N} \bm u_{n,1} \rVert_{2}^2\\
\leq & \frac{2}{\varepsilon^2}\lVert \nabla_{N}(\bm u^n)^3 \rVert_{2}^2 + \frac{\varepsilon^2}{2}\lVert  \left(\phi_{1}(\tau L_{\kappa})\right)^{\frac{1}{2}} \nabla_{N}\Delta_{N} \bm u_{n,1}  \rVert_{2}^2 \\
\leq & \frac{2C}{\varepsilon^2} \lVert \bm u^n \rVert_{\infty}^4 \lVert \nabla_{N}\bm u^n \rVert_{2}^2 + \frac{\varepsilon^2}{2} \lVert \left(\phi_{1}(\tau L_{\kappa})\right)^{\frac{1}{2}} \nabla_{N}\Delta_{N} \bm u_{n,1}  \rVert_{2}^2 \\
\leq & \frac{2C}{\varepsilon^2} \tilde{C}_n^4 \lVert \nabla_{N} \bm u^n \rVert_{2}^2 + \frac{\varepsilon^2}{2} \lVert \left(\phi_{1}(\tau L_{\kappa})\right)^{\frac{1}{2}} \nabla_{N}\Delta_{N} \bm u_{n,1}  \rVert_{2}^2.
\end{aligned}
\end{equation}
The second term on the right-hand side of \eqref{nonl_inn} can be treated similarly. It holds that
\begin{equation}\label{eqn:rough1_h1_8}
  \begin{aligned}
  & 2 \kappa  \left< \phi_{1}(\tau L_{\kappa}) \Delta_{N} \bm u^n , \Delta_{N}\bm u_{n,1} \right> \\
  = & \kappa (\lVert \left(\phi_{1}(\tau L_{\kappa})\right)^{\frac{1}{2}} \Delta_{N} \bm u^n  \rVert_{2}^2 + \lVert \left(\phi_{1}(\tau L_{\kappa})\right)^{\frac{1}{2}}\Delta_N \bm u_{n,1} \rVert^2_{2} - \lVert \left(\phi_{1}(\tau L_{\kappa})\right)^{\frac{1}{2}} \Delta_{N} (\bm u_{n,1} - \bm u_{n}) \rVert_2^2 ).
  \end{aligned}
\end{equation}
Furthermore, by Young's inequality and the operator estimate \eqref{eqn:phi_esitmate}, the remaining part of \eqref{nonl_inn} can be bounded by
\begin{equation}\label{eqn:rough1_h1_9}
  2\left< \phi_{1}(\tau L_{\kappa}) \Delta_{N} \bm u^n , \Delta_{N} \bm u_{n,1} \right> \leq \frac{2}{\varepsilon^2} \lVert \nabla_{N} \bm u^n \rVert_{2}^2 + \frac{\varepsilon^2}{2} \lVert  \left(\phi_{1}(\tau L_{\kappa})\right)^{\frac{1}{2}} \nabla_{N} \Delta_{N} \bm u_{n,1} \rVert_{2}^2.
\end{equation}
Substituting \eqref{eqn:rough1_h1_7}--\eqref{eqn:rough1_h1_9} into \eqref{nonl_inn}, and combining the result with \eqref{ineq_G1} and \eqref{eqn:rough1_h1_5}, we obtain
\[
  \begin{aligned}
\lVert \nabla_{N} \bm u_{n,1} \rVert_{2}^2 - \lVert \nabla_{N}\bm u^n \rVert_{2}^2 +& \varepsilon^2\tau \left\lVert \left(\phi_{1}(\tau L_{\kappa})\right)^{\frac{1}{2}} \nabla_N \Delta_N \bm u^n \right\rVert_{2}^2 \\
+ &  \kappa \tau \lVert (\phi_{1}(\tau L_\kappa))^{\frac{1}{2}}\Delta_{N} (\bm u_{n,1} - \bm u_{n}) \rVert_2^2  \leq  \frac{2\tau}{\varepsilon^{2}}(C \tilde{C}_n^4 +1 ) \lVert \nabla_{N}\bm u^n \rVert_{2}^2.
\end{aligned}
\]
Consequently, we have
\[
    \left\lVert \nabla_{N} \bm u_{n,1} \right\rVert_{2}^2 + \kappa \tau \left\lVert \phi_{1}^{\frac{1}{2}}(\tau L_\kappa)\Delta_{N} (\bm u_{n,1} - \bm u_{n}) \right\rVert_2^2 \leq \left(1+ \frac{2\tau}{\varepsilon^{2}} (C \tilde{C}_n^4 +1 )\right) \left\lVert \nabla_{N}\bm u^n \right\rVert_{2}^2.
\]
Under the time step restriction
\begin{equation}\label{eqn:stepConst_1}
  \tau \leq \tau_1 := \varepsilon^2 (2\del{ }C \tilde{C}_n^4 + 2)^{-1},
\end{equation}
we obtain the preliminary bound
\[\left\Vert \nabla_N \bm u_{n,1} \right\Vert_2^2 \leq 2 \left\Vert \nabla_N \bm u^n \right\Vert_2^2,\]
which implies
\begin{equation}\label{eqn:stage1_H1coarseEst}
\left\Vert \nabla_N \bm u_{n,1} \right\Vert_2 \leq \sqrt{2} \left\Vert \nabla_N \bm u^n \right\Vert_2 \leq \sqrt{2}\dot{C}_n,
\end{equation}
where the last inequality follows from the a priori estimate \eqref{eqn:apriori_H1}.
Moreover, the artificial dissipation term satisfies the same bound:
\begin{equation}\label{eqn:stage1_articifialEst}
  \kappa \tau \lVert (\phi_{1}(\tau L_\kappa))^{\frac{1}{2}}\Delta_{N} (\bm u_{n,1} - \bm u^{n}) \rVert_2^2 \leq 2 \left\Vert \nabla_N \bm u^n \right\Vert_2^2,
\end{equation}
which will facilitate the analysis in the next stage.

\subsection{Preliminary estimate for \texorpdfstring{$\Vert\Delta_N \bm u_{n,1}\Vert_2$}{}}\label{Xsec9-3.3}
By taking the discrete $\ell^2$ inner product of both sides of \eqref{eqn:rfd_1} with $2\Delta_N^2 \bm u_{n,1}$ yields
\begin{equation}\label{eqn:rough1_h2_1}
  \left\langle \bm u_{n,1} - \bm u^n , 2\Delta_N^2 \bm u_{n,1} \right\rangle  + 2\tau \left\langle \phi_1(\tau L_\kappa) L_\kappa \bm u^n, \Delta_N^2 \bm u_{n,1} \right\rangle = 2 \tau \left\langle \phi_1(\tau L_\kappa) N_\kappa(\bm u^n), \Delta_N^2 \bm u_{n,1} \right\rangle.
\end{equation}
Following the derivation in \eqref{eqn:rough1_h1_2}--\eqref{eqn:rough1_h1_3}, we can rewrite the linear terms on the left-hand side as
\[
  \left\langle \bm u_{n,1} - \bm u^n , 2\Delta_N^2 \bm u_{n,1} \right\rangle = \left\Vert \Delta_N \bm u_{n,1} \right\Vert_2^2 - \left\Vert \Delta_N \bm u^n \right\Vert_2^2 + \left\Vert \Delta_N (\bm u_{n,1} - \bm u^n) \right\Vert_2^2
\]
and
\[
  \begin{aligned}
  2 \left\langle \phi_1(\tau L_\kappa) L_\kappa \bm u^n, \Delta_N^2 \bm u_{n,1}  \right\rangle =& 2\langle G_N^{(2)} \bm u^n , G_N^{(2)} \bm u_{n,1} \rangle \\ =& \lVert G_{N}^{(2)}\bm u^n \rVert_{2}^2 +  \lVert G_{N}^{(2)} \bm u_{n,1} \rVert_{2}^2 - \lVert G_{N}^{(2)} (\bm u^n - \bm u_{n,1}) \rVert_{2}^2.
  \end{aligned}
\]
Additionally, by \eqref{eqn:diffusion_estimate}, it holds that
\[
  \left\Vert \Delta_N (\bm u_{n,1} - \bm u^n) \right\Vert_2^2 - \tau  \lVert G_{N}^{(2)} (\bm u^n - \bm u_{n,1}) \rVert_{2}^2 \geq 0.
\]
Substituting these identities and the above inequality into \eqref{eqn:rough1_h2_1} yields
\begin{equation}\label{eqn:rough1_h2_2}
\left\Vert \Delta_N \bm u_{n,1} \right\Vert_2^2 - \left\Vert \Delta_N \bm u^n \right\Vert_2^2 + \tau (\lVert G_N^{(2)} \bm u^n \rVert_{2}^2 +  \lVert G_{N}^{(2)} \bm u_{n,1} \rVert_{2}^2) \leq 2 \tau \left\langle \phi_1(\tau L_\kappa) N_\kappa(\bm u^n), \Delta_N^2 \bm u_{n,1} \right\rangle.
\end{equation}
The two terms involving $G_N^{(2)}$ can be expanded as
\[
\begin{aligned}
    \lVert G_N^{(2)} \bm u^n \rVert_{2}^2 +  \lVert G_{N}^{(2)} \bm u_{n,1} \rVert_{2}^2 =& \varepsilon^2 (\Vert \left( \phi_1(\tau L_\kappa)\right)^{\frac{1}{2}} \Delta_N^2 \bm u^{n} \Vert_2^2 + \Vert \left( \phi_1(\tau L_\kappa)\right)^{\frac{1}{2}} \Delta_N^2 \bm u_{n,1} \Vert_2^2)\\
    &+ \kappa (\| (\phi_1(\tau L_\kappa))^{\frac12} \Delta_N \nabla_N \bm u^n \|_2^2 + \| (\phi_1(\tau L_\kappa))^{\frac12} \Delta_N \nabla_N \bm u_{n,1} \|_2^2).
\end{aligned}
\]
Similarly, the nonlinear term on the right-hand side of \eqref{eqn:rough1_h2_2} can be expressed as
\begin{equation}\label{eqn:rough1_h2_3}
  \begin{aligned}
  &2\left\langle \phi_1(\tau L_\kappa) N_\kappa(\bm u^n), \Delta_N^2 \bm u_{n,1} \right\rangle  \\
 =&2 \left< \phi_{1}(\tau L_{\kappa})\Delta_{N} (\bm u^n)^3, \Delta_{N}^2 \bm u_{n,1} \right> - 2(\kappa + 1) \left< \phi_{1}(\tau L_{\kappa}) \Delta_{N} \bm u^n , \Delta_{N}^2 \bm u_{n,1} \right>.
  \end{aligned}
\end{equation}
Applying estimates \eqref{eqn:nonlineEst3} and \eqref{eqn:phi_esitmate} to the right-hand side of \eqref{eqn:rough1_h2_3}, we obtain
\begin{equation}\label{eqn:rough1_h2_31}
\begin{aligned}
2\left< \phi_{1}(\tau L_{\kappa}) \Delta_N (\bm u^n)^3, \Delta_N^2 \bm u_{n,1} \right> =& 2\left< \left(\phi_{1}(\tau L_{\kappa})\right)^{\frac{1}{2}} \Delta_N (\bm u^n)^3,\left(\phi_{1}(\tau L_{\kappa})\right)^{\frac{1}{2}}\Delta_N^2 \bm u_{n,1} \right> \\
\leq & \frac{2}{\varepsilon^2}\left\lVert \Delta_N (\bm u^n)^3 \right\rVert_{2}^2 + \frac{\varepsilon^2}{2}\lVert \left(\phi_{1}(\tau L_{\kappa})\right)^{\frac{1}{2}}\Delta_N^2 \bm u_{n,1} \rVert_{2}^2 \\
\leq & \frac{2C}{\varepsilon^2}\left\lVert \bm u^n \right\rVert_{\infty}^4\left\lVert \Delta_N \bm u^n \right\rVert_{2}^2 + \frac{\varepsilon^2}{2}\lVert \left(\phi_{1}(\tau L_{\kappa})\right)^{\frac{1}{2}}\Delta_N^2 \bm u_{n,1} \rVert_{2}^2\\
\leq & \frac{2C}{\varepsilon^2}\tilde{C}_n^4\left\lVert \Delta_N \bm u^n \right\rVert_{2}^2 + \frac{\varepsilon^2}{2}\lVert \left(\phi_{1}(\tau L_{\kappa})\right)^{\frac{1}{2}}\Delta_N^2 \bm u_{n,1}  \rVert_{2}^2
\end{aligned}
\end{equation}
and
\begin{equation}\label{eqn:rough1_h2_32}
\begin{aligned}
-2(\kappa + 1) \left< \phi_{1}(\tau L_{\kappa}) \Delta_N \bm u^{n}, \Delta_N^2 \bm u_{n,1} \right>
\leq & \kappa ( \lVert  \left(\phi_{1}(\tau L_{\kappa})\right)^{\frac{1}{2}} \nabla_{N} \Delta_N \bm u^n  \rVert_{2}^2 + \lVert  \left(\phi_{1}(\tau L_{\kappa})\right)^{\frac{1}{2}} \nabla_N \Delta_N \bm u_{n,1} \rVert_{2}^2  )   \\
&+ \frac{2}{\varepsilon^2} \left\lVert \Delta_N \bm u^n  \right\rVert_{2}^2 + \frac{\varepsilon^2}{2} \lVert \phi_{1}^{ \frac{1}{2} } (\tau L_{\kappa}) \Delta_N^2 \bm u_{n,1}  \rVert_{2}^2.
\end{aligned}
\end{equation}
Combining and rearranging the estimates  {\eqref{eqn:rough1_h2_2}--\eqref{eqn:rough1_h2_32}} yields
\[
\lVert \Delta_{N} \bm u_{n,1} \rVert^2 - \lVert \Delta_{N} \bm u^n \rVert^2 + \varepsilon^2 \tau \lVert \left(\phi_{1}(\tau L_{\kappa})\right)^{\frac{1}{2}} \Delta_N^2 \bm u^n \rVert_{2}^2  \leq  2\varepsilon^{-2}\tau \left( C\tilde{C}_n^4 + 1 \right)\left\lVert \Delta_N \bm u^n \right\rVert_{2}^2.
\]
This leads to the bound
\[
\lVert \Delta_{N} \bm u_{n,1} \rVert_2^2 \leq  \left(1 + 2\varepsilon^{-2}\tau \left( C {\tilde{C}_n^4}+ 1 \right)\right)\left\lVert \Delta_N \bm u^n  \right\rVert_{2}^2.
\]
Under the time step constraint \eqref{eqn:stepConst_1}, we obtain the preliminary estimate
\begin{equation}\label{eqn:stage1_H2coarseEst}
\lVert \Delta_{N} \bm u_{n,1} \rVert_2^2  \leq  2 \lVert \Delta_{N} \bm u^n \rVert_2^2,\ \text{so that}\ \lVert \Delta_{N} \bm u_{n,1} \rVert_2  \leq  \sqrt{2} \lVert \Delta_{N} \bm u^n \rVert_2 \leq \sqrt{2} \ddot{C}_n,
\end{equation}
where the last inequality follows from the a priori estimate \eqref{eqn:apriori_H2}.

Combining \eqref{eqn:stage1_H1coarseEst} and \eqref{eqn:stage1_H2coarseEst} and using the discrete Sobolev inequality \eqref{eqn:3dSob} for the discrete maximum norm, we further obtain
\[
\begin{aligned}
  \left\Vert \bm u_{n,1} \right\Vert_\infty \leq & \tilde{C}'_{n,1} := C\left\Vert \nabla_N \bm u_{n,1} \right\Vert^{\frac{1}{2}} \left\Vert \Delta_N \bm u_{n,1} \right\Vert^{\frac{1}{2}}.
\end{aligned}
\]
This, in turn, allows us to establish a rough bound for the intermediate solution $\bm u_{n,1}$ in  $\ell^\infty$ norm. Furthermore, by \linkref[Theorem]{\ref{thm:p_stab}} and \linkref[Lemma]{\ref{lem:energy_est}}, there exists a stabilization parameter $\kappa$ that guarantees the energy stability of the first stage in the numerical scheme, i.e., $E_N(\bm u_{n,1}) \leq E_{N}(\bm u^n)$. Leveraging this, we derive a refined estimate of the $H^1$ semi-norm and the $\ell^{\infty}$ norm of the intermediate solution $\bm u_{n,1}$:
\begin{equation}\label{eqn:stage1_H1refineEst}
        \Vert \nabla_N \bm u_{n,1} \Vert_2 \leq \dot{C}_n,\quad \Vert \bm u_{n,1} \Vert_\infty \leq \tilde{C}_{n,1} \leq \sqrt[4]{2}C\dot{C}_n^{\frac{1}{2}} \ddot{C}_n^{\frac{1}{2}}.
\end{equation}

\subsection{Preliminary estimate for \texorpdfstring{$\Vert\nabla_N \bm u^{n+1}\Vert_2$}{}}\label{Xsec10-3.4}
Next, we analyze the numerical solution $\bm u^{n+1}$ generated by the second stage of the ERK2 scheme \eqref{eqn:etd_s1}--\eqref{eqn:etd_s2}. Taking the discrete $\ell^2$ inner product of \eqref{eqn:rfd_2} with $-2\Delta_N \bm u^{n+1}$ yields
\[
  \begin{aligned}
    &\left\langle \bm u^{n+1} - \bm u^n, -2\Delta_N \bm u^{n+1} \right\rangle - 2\tau \left\langle \phi_1(\tau L_\kappa) L_\kappa \bm u^n, \Delta_N \bm u^{n+1} \right\rangle\\ =& -\tau \left\langle \phi_1(\tau L_\kappa) N_\kappa\left(\bm u^n\right), \Delta_N \bm u^{n+1} \right\rangle - \tau \left\langle \phi_1(\tau L_\kappa) N_\kappa(\bm u_{n,1}), \Delta_N \bm u^{n+1} \right\rangle.
  \end{aligned}
\]
The temporal difference and diffusion terms on the left-hand side of the above identity can be estimated as in \eqref{eqn:rough1_h1_2}--\eqref{eqn:rough1_h1_4}, leading to
\begin{equation}\label{eqn:stage2_h1_1}
  \begin{aligned}
    &\left\Vert \nabla_N \bm u^{n+1} \right\Vert_2^2 - \left\Vert \nabla_N \bm u^n \right\Vert_2^2 + \tau (\Vert G_N^{(1)} \bm u^n \Vert_2^2 + \Vert G_N^{(1)} \bm u^{n+1} \Vert_2^2 )\\
    \leq & - \tau \left\langle \phi_1(\tau L_\kappa) N_\kappa(\bm u^n), \Delta_N \bm u^{n+1} \right\rangle - \tau \left\langle \phi_1(\tau L_\kappa) N_\kappa(\bm u_{n,1}), \Delta_N \bm u^{n+1} \right\rangle,
  \end{aligned}
\end{equation}
where {the terms involving the operator $G_N^{(1)}$ from \eqref{eqn:linear_operators} can be expressed as}
\[
\begin{aligned}
    \Vert G_N^{(1)} \bm u^n \Vert_2^2 + \Vert G_N^{(1)} \bm u^{n+1} \Vert_2^2
    =& \varepsilon^2 (\Vert \left( \phi_1(\tau L_\kappa)\right)^{\frac{1}{2}} \Delta_N \nabla_N \bm u^n \Vert_2^2 + \Vert \left( \phi_1(\tau L_\kappa)\right)^{\frac{1}{2}} \Delta_N \nabla_N \bm u^{n+1} \Vert_2^2) \\ &+ \kappa (\Vert \left( \phi_1(\tau L_\kappa)\right)^{\frac{1}{2}} \Delta_N \bm u^n \Vert_2^2  + \Vert \left( \phi_1(\tau L_\kappa)\right)^{\frac{1}{2}} \Delta_N \bm u^{n+1} \Vert_2^2).
\end{aligned}
\]
The first term on the right-hand side of \eqref{eqn:stage2_h1_1} can be decomposed using the same approach as in \eqref{eqn:rough1_h1_7}--\eqref{eqn:rough1_h1_9}. For brevity, we omit the details and present the result directly:
\begin{equation}\label{eqn:stage2_h1_2}
  \begin{aligned}
   &- \left\langle \phi_1(\tau L_\kappa) N_\kappa(\bm u^n), \Delta_N \bm u^{n+1} \right\rangle \\
  =& - \left\langle \phi_1(\tau L_\kappa) \Delta_N(\bm u^n)^3, \Delta_N \bm u^{n+1} \right\rangle + (\kappa + 1)\left\langle \phi_1(\tau L_\kappa) \Delta_N \bm u^n, \Delta_N \bm u^{n+1} \right\rangle \\
  \leq & \frac{\varepsilon^2}{2}  \lVert \left(\phi_{1}(\tau L_{\kappa})\right)^{\frac{1}{2}} \nabla_{N}\Delta_{N} \bm u^{n+1}  \rVert_{2}^2 + \frac{\kappa}{2} (\lVert \left(\phi_{1}(\tau L_{\kappa})\right)^{\frac{1}{2}} \Delta_{N} \bm u^n  \rVert_{2}^2 + \lVert \left(\phi_{1}(\tau L_{\kappa})\right)^{\frac{1}{2}}\Delta \bm u^{n+1} \rVert^2_{2} ) \\
  & + \frac{1 + C\tilde{C}_n^4}{\varepsilon^2} \lVert \nabla_{N} \bm u^n \rVert_{2}^2  - \frac{\kappa}{2} \lVert \left(\phi_{1}(\tau L_{\kappa})\right)^{\frac{1}{2}}\Delta_N( \bm u^{n+1} - \bm u^n )\rVert^2_{2}.
  \end{aligned}
\end{equation}
The remaining nonlinear inner product is treated analogously, with an additional artificial dissipation term introduced via the triangular inequality:
\begin{equation}\label{eqn:stage2_h1_3}
  \begin{aligned}
   & - \left\langle \phi_1(\tau L_\kappa) N_\kappa(\bm u_{n,1}), \Delta_N \bm u^{n+1} \right\rangle \\
  =& - \left\langle \phi_1(\tau L_\kappa) \Delta_N(\bm u_{n,1})^3, \Delta_N \bm u^{n+1} \right\rangle + (\kappa + 1)\left\langle \phi_1(\tau L_\kappa) \Delta_N \bm u_{n,1}, \Delta_N \bm u^{n+1} \right\rangle \\
  \leq & \frac{\varepsilon^2}{2}  \lVert \left(\phi_{1}(\tau L_{\kappa})\right)^{\frac{1}{2}} \nabla_{N}\Delta_{N} \bm u^{n+1}  \rVert_{2}^2 + \frac{\kappa}{2} (\lVert \left(\phi_{1}(\tau L_{\kappa})\right)^{\frac{1}{2}} \Delta_{N} \bm u^{n} \rVert_{2}^2 + \lVert \left(\phi_{1}(\tau L_{\kappa})\right)^{\frac{1}{2}}\Delta_N \bm u^{n+1} \rVert^2_{2} ) \\
  & + \frac{1 + C\tilde{C}_{n,1}^4}{\varepsilon^2} \lVert \nabla_{N} \bm u_{n,1} \rVert_{2}^2 - \frac{\kappa}{2} \lVert \left(\phi_{1}(\tau L_{\kappa})\right)^{\frac{1}{2}}\Delta_N( \bm u^{n+1} - \bm u^n )\rVert^2_{2} + \kappa \left\langle \phi_1(\tau L_\kappa) \Delta_N (\bm u_{n,1} -  \bm u^{n}) , \Delta_N \bm u^{n+1} \right\rangle.
  \end{aligned}
\end{equation}
The artificial dissipation term can be estimated as
\begin{equation}
\label{eqn:stage2_h1_3disst}
\begin{aligned}
     &\kappa \left\langle \phi_1(\tau L_\kappa) \Delta_N (\bm u_{n,1} -  \bm u^{n}) , \Delta_N \bm u^{n+1} \right\rangle\\
     \leq& \frac{\kappa}{2} \Vert (\phi_1(\tau L_\kappa))^{\frac{1}{2}} \Delta_N (\bm u_{n,1} - \bm u^{n}) \Vert_2^2 + \frac{\kappa}{2} \Vert (\phi_1(\tau L_\kappa))^{\frac{1}{2}} \Delta_N \bm u^{n+1} \Vert_2^2\\
     \leq & \frac{\kappa}{2} \Vert (\phi_1(\tau L_\kappa))^{\frac{1}{2}} \Delta_N (\bm u_{n,1} - \bm u^{n}) \Vert^2 + \kappa \Vert (\phi_1(\tau L_\kappa))^{\frac{1}{2}} \Delta_N \bm u^{n} \Vert^2 + \kappa \Vert (\phi_1(\tau L_\kappa))^{\frac{1}{2}} \Delta_N (\bm u^{n+1} -  \bm u^{n}) \Vert^2,
\end{aligned}
\end{equation}
where the first term on the right-hand side can be bounded using \eqref{eqn:stage1_articifialEst}. For the second term, we apply   {Young's} inequality to obtain
\begin{equation}
\label{eqn:stage2_h1_3kap}
    \kappa\Vert (\phi_1(\tau L_\kappa))^{\frac{1}{2}} \Delta_N \bm u^{n} \Vert_2^2 \leq \varepsilon^2 \Vert (\phi_1(\tau L_\kappa))^{\frac{1}{2}} \Delta_N \nabla_N \bm u^{n} \Vert_2^2 + \frac{\kappa^2}{4\varepsilon^2}\Vert \nabla_N \bm u^n \Vert_2^2.
\end{equation}
Therefore, by {combining \eqref{eqn:stage2_h1_1}--\eqref{eqn:stage2_h1_3kap} with }   the  preliminary estimate \eqref{eqn:stage1_articifialEst} from the previous stage, we have
\[
  \begin{aligned}
  \lVert \nabla_{N} \bm u^{n+1} \rVert_{2}^2 - \lVert \nabla_{N} \bm u^n \rVert_{2}^2
  \leq \left(1 + (1 + C\tilde{C}_n^4 + \frac{\kappa^2}{4})\varepsilon^{-2} \tau \right) \lVert \nabla_{N} \bm u^n \rVert_{2}^2 + \left(1 + C \tilde{C}_{n,1}^4\right) \varepsilon^{-2} \tau \lVert \nabla_{N} \bm u_{n,1} \rVert_{2}^2.
\end{aligned}
\]
From the inequality above, it follows that
\[
\begin{aligned}
 \lVert \nabla_{N} \bm u^{n+1} \rVert_{2}^2 \leq & (2 +  \varepsilon^{-2}\tau(C \tilde{C}_n^4 + \frac{\kappa^2}{4} + 1)) \lVert \nabla_{N} \bm u^n \rVert_{2}^2 + \varepsilon^{-2}\tau (C \tilde{C}_{n}^4 + 1) \lVert \nabla_{N} \bm u_{n,1} \rVert_{2}^2.
  \end{aligned}
\]
By combining this with  the estimate \eqref{eqn:stage1_H1refineEst} and assuming
\begin{equation}\label{eqn:stepConst_2}
  \tau \leq \tau_2 := \varepsilon^2 \left(C \tilde{C}_n^4 + \frac{\kappa^2}{4} + 1\right)^{-1},
\end{equation}
we can derive
\[
  \begin{aligned}
    \lVert \nabla_{N} \bm u^{n+1} \rVert_{2}^2 \leq 3\lVert \nabla_{N} \bm u^n \rVert_{2}^2 + \lVert \nabla_{N} \bm u_{n,1} \rVert_{2}^2 \leq 4\dot{C}_n^2, \text{ so that } \lVert \nabla_{N} \bm u^{n+1} \rVert_{2} \leq 2 \dot{C}_n.
  \end{aligned}
\]

\subsection{Uniform-in-time estimate for \texorpdfstring{$\Vert\Delta_N u^{n+1}\Vert_2$}{}}\label{Xsec11-3.5}
It is observed that the a priori assumption given in \eqref{eqn:apriori_H1} is uniform-in-time and can be recovered via the energy stability analysis at the next time step. By contrast, the theoretical justification for the a priori assumption in \eqref{eqn:apriori_H2} proves more challenging, as the global-in-time constants $\widehat{C}_{\varepsilon,0}$ and $\widehat{C}_{\varepsilon,1}$ have not yet been derived. In this part, we aim to establish an $H^2$ estimate for the numerical solution $\bm u^{n+1}$ at the next time step. Specifically, a relatively rough $H^2$ estimate for the intermediate solution $\bm u_{n,1}$ has been derived in \eqref{eqn:stage1_H2coarseEst}. However, due to the nonlinear structure of the CH equation, the same analytical technique employed for the formulation in \eqref{eqn:rfd_2} is not well-suited for recovering the a priori estimate given at the previous time step \eqref{eqn:apriori_H2}.

To obtain a more refined estimate, we first apply the inverse operator $(\phi{_1}(\tau L_\kappa))^{-1}$ to the numerical solution $\bm{u}^{n+1}$. This allows us to rewrite \eqref{eqn:rfd_2} in the alternative form:
\begin{equation}\label{eqn:stage2_rfd}
  (\phi_1(\tau L_\kappa) )^{-1} (\bm u^{n+1} - \bm u^n) +  \tau L_\kappa \bm u^n = \frac{1}{2}\tau (N_\kappa(\bm u^n) + N_\kappa(\bm u_{n,1})).
\end{equation}
Taking the discrete $\ell^2$ inner product of both sides of \eqref{eqn:stage2_rfd} with $2\Delta_N^2 u^{n+1}$, we derive
\[
  \begin{aligned}
   \left\langle (\phi_1(\tau L_\kappa))^{-1} (\bm u^{n+1} - \bm u^n), 2\Delta_N^2 \bm u^{n+1} \right\rangle + 2\tau \left\langle L_\kappa \bm u^n, \Delta_N^2 \bm u^{n+1} \right\rangle
  = \frac{1}{2}\tau \left\langle  (N_\kappa(\bm u^n) + N_\kappa(\bm u_{n,1})), 2\Delta_N^2 \bm u^{n+1} \right\rangle.
  \end{aligned}
\]
The linear terms on the left-hand side can be estimated using an approach similar to that in \eqref{eqn:rough1_h1_2}--\eqref{eqn:rough1_h1_4}, yielding
\begin{equation}\label{eqn:stage2_h2_1}
  \begin{aligned}
  &\left\langle (\phi_1 (\tau L_\kappa) )^{-1}(\bm u^{n+1} - \bm u^n), 2\Delta_N^2 \bm u^{n+1} \right\rangle \\
  =& \Vert (\phi_1 (\tau L_\kappa) )^{-\frac{1}{2}} \Delta_N \bm u^{n+1} \Vert _2^2 - \Vert (\phi_1(\tau L_\kappa))^{-\frac{1}{2}} \Delta_N \bm u^n \Vert_2^2 + \Vert (\phi_1(\tau L_\kappa))^{-\frac{1}{2}} \Delta_N (\bm u^{n+1} - \bm u^n) \Vert_2^2
  \end{aligned}
\end{equation}
and
\begin{equation}\label{eqn:stage2_h2_2}
  2\left\langle L_\kappa \bm u^n, \Delta_N^2 \bm u^{n+1} \right\rangle = \Vert \tilde{G}_N^{(2)}\bm u^n \Vert_2^2 + \Vert \tilde{G}_N^{(2)} \bm u^{n+1} \Vert_2^2 - \Vert \tilde{G}_N^{(2)} (\bm u^{n+1} - \bm u^{n}) \Vert_2^2.
\end{equation}
Next, applying the definition of $\tilde{G}_N^{(2)}$ and the result from \eqref{eqn:diffusion_estimate} gives
\begin{align}\label{eqn:stage2_h2_3}
&   \Vert \tilde{G}_N^{(2)} \bm u^n \Vert_2^2 + \Vert \tilde{G}_N^{(2)}\bm u^{n+1} \Vert_2^2 = \varepsilon^2 (\left\Vert \Delta^2_N \bm u^n \right\Vert_2^2  + \left\Vert \Delta^2_N \bm u^{n+1} \right\Vert_2^2 ) + \kappa (\left\Vert \Delta_N \nabla_N \bm u^n \right\Vert_2^2 +\left\Vert \Delta_N \nabla_N \bm u^{n+1} \right\Vert_2^2 ),\\
  \label{eqn:stage2_h2_4}
&   \lVert (\phi_1(\tau L_{\kappa}))^{-\frac{1}{2}} \Delta_N (\bm u^{n+1} - \bm u^n) \rVert_{2}^2 - \tau \lVert \tilde{G}_{N}^{(2)} (\bm u^{n+1} -\bm u^n) \rVert_{2}^2 \geq 0.
\end{align}
Therefore, by combining \eqref{eqn:stage2_h2_1}, \eqref{eqn:stage2_h2_2} and\eqref{eqn:stage2_h2_4}, we obtain
\begin{equation}\label{eqn:stage2_h2_5}
  \begin{aligned}
  &\Vert (\phi_1(\tau L_\kappa))^{-\frac{1}{2}} \Delta_N \bm u^{n+1} \Vert _2^2 - \Vert (\phi_1(\tau L_\kappa))^{-\frac{1}{2}} \Delta_N \bm u^n \Vert_2^2 + \tau \Vert \tilde{G}_N^{(2)}\bm u^n \Vert_2^2 + \tau \Vert \tilde{G}_N^{(2)}\bm u^{n+1} \Vert_2^2 \\
  \leq & \tau \langle  (N_\kappa(\bm u^n) + N_\kappa(\bm u_{n,1})), \Delta_N^2 \bm u^{n+1} \rangle \\
  = & \tau \left\langle \Delta_N (\bm u^n)^3 , \Delta_N^2 \bm u^{n+1} \right\rangle + \tau\left\langle \Delta_N (\bm u_{n,1})^3 , \Delta_N^2 \bm u^{n+1} \right\rangle \\
  & -\tau(\kappa + 1)\left\langle \Delta_N \bm u^n, \Delta_N^2 \bm u^{n+1} \right\rangle - \tau(\kappa + 1)\left\langle \Delta_N \bm u_{n,1}, \Delta_N^2 \bm u^{n+1} \right\rangle.
  \end{aligned}
\end{equation}
Regarding the first two nonlinear inner product terms on the right-hand side of \eqref{eqn:stage2_h2_5}, {the nonlinear estimate \eqref{eqn:nonlineEst3}, together with the discrete Sobolev inequality \eqref{eqn:3dSob}, is applied in a manner similar to the derivation of \eqref{eqn:rough1_h1_7}. This yields}
\begin{equation}\label{eqn:stage2_h2_6}
\scalebox{0.92}{
$\begin{aligned}
\left\langle \Delta_N (\bm u^{n})^3 , \Delta_N^2 \bm u^{n+1} \right\rangle
&\leq  \frac{2}{\varepsilon^{2}} \left\Vert \Delta_N (\bm u^n)^3 \right\Vert_2^2 + \frac{\varepsilon^{2}}{8} \left\Vert \Delta_N^2 \bm u^{n+1} \right\Vert_2^2 \\
&\leq \frac{2C}{\varepsilon^{2}} \left\Vert \nabla_N \bm u^n \right\Vert_2^2 \left\Vert \Delta_N \bm u^n \right\Vert_2^4 + \frac{\varepsilon^{2}}{8} \left\Vert \Delta_N^2 \bm u^{n+1} \right\Vert_2^2\\
&\leq  \frac{2C}{\varepsilon^{2}} \dot{C}_n^2 \left\Vert \Delta_N \bm u^n \right\Vert_2^4 + \frac{\varepsilon^{2}}{8} \left\Vert \Delta_N^2 \bm u^{n+1} \right\Vert_2^2
\end{aligned}$
}
\end{equation}
and
\begin{equation}\label{eqn:stage2_h2_7}
\resizebox{0.92\textwidth}{!}{$
  \begin{aligned}
  \left\langle \Delta_N (\bm u_{n,1})^3, \Delta_N^2 \bm u^{n+1} \right\rangle \leq & \frac{2}{\varepsilon^{2}}\lVert \Delta_N (\bm u_{n,1})^3 \rVert_{2}^2 + \frac{\varepsilon^2}{8} \lVert \Delta_N^2 \bm u^{n+1} \rVert_{2}^2
  \leq \frac{8C}{\varepsilon^{2}} \lVert \nabla_N \bm u_{n,1} \rVert_{2}^2  \lVert \Delta_N \bm u^{n} \rVert_{2}^4 + \frac{\varepsilon^2}{8} \lVert \Delta_N^2 \bm u^{n+1} \rVert_{2}^2\\
  \leq & \frac{8C}{\varepsilon^{2}} \dot{C}_n^2  \lVert \Delta_N \bm u^{n} \rVert_{2}^4 + \frac{\varepsilon^2}{8} \lVert \Delta_N^2 \bm u^{n+1} \rVert_{2}^2.
  \end{aligned}$}
\end{equation}
For the third term in \eqref{eqn:stage2_h2_5}, a similar bound can be established by means of the summation by parts formula and Young's inequality:
\begin{equation}\label{eqn:stage2_h2_8}
  -(\kappa +1) \left\langle \Delta \bm u^n, \Delta_N^2 \bm u^{n+1} \right\rangle \leq  \frac{2}{\varepsilon^{2}} \left\Vert \Delta_N \bm u^n \right\Vert_2^2 + \frac{\varepsilon^2}{8} \left\Vert \Delta^2_N \bm u^{n+1} \right\Vert_2^2 + \frac{\kappa}{2}(\left\Vert \nabla_N \Delta_N \bm u^n \right\Vert_2^2 + \left\Vert \nabla_N \Delta_N \bm u^{n+1} \right\Vert_2^2).
\end{equation}
For the last term on the right-hand side of \eqref{eqn:stage2_h2_5}, when combined with the preliminary estimate associated with $\|\Delta_N \bm u_{n,1}\|_2$, it can be analyzed as follows:
\begin{equation}\label{eqn:stage2_h2_9}
  \begin{aligned}
  - (\kappa + 1) \left< \Delta_N \bm u_{n,1}, \Delta_N^2 \bm u^{n+1} \right> \leq & \frac{2(\kappa^2+1)}{\varepsilon^{2}} \lVert \Delta_N \bm u_{n,1} \rVert_{2}^2 + \frac{\varepsilon^2}{8} \lVert \Delta_N^2 \bm u^{n+1} \rVert_{2}^2 \\
  \leq & \frac{4(\kappa^2+1)}{\varepsilon^{2}} \lVert \Delta_N \bm u_{n} \rVert_{2}^2 + \frac{\varepsilon^2}{8} \lVert \Delta_N^2 \bm u^{n+1} \rVert_{2}^2.
  \end{aligned}
\end{equation}
Subsequently, combining \eqref{eqn:stage2_h2_6}--\eqref{eqn:stage2_h2_9} with \eqref{eqn:stage2_h2_5} yields
\begin{equation}\label{eqn:stage2_h2_10}
  \begin{aligned}
    &\Vert (\phi_1(\tau L_\kappa))^{-\frac{1}{2}}  \Delta_N \bm u^{n+1} \Vert_2^2 - \Vert (\phi_1(\tau L_\kappa))^{-\frac{1}{2}} \Delta_N \bm u^n \Vert_2^2 + \frac{\varepsilon^2}{2} \tau \Vert \Delta_N^2 \bm u^{n+1} \Vert_2^2 + \frac{\kappa}{2}\left\Vert \nabla_N \Delta_N \bm u^{n+1} \right\Vert_2^2 \\
   & + \varepsilon^2 \tau \Vert \Delta_N^2 \bm u^n \Vert_2^2 + \frac{\kappa}{2}\left\Vert \nabla_N \Delta_N \bm u^n \right\Vert_2^2  \leq  {2(\kappa^2+3)} \varepsilon^{-2}\tau \left\Vert \Delta_N \bm u^n \right\Vert_2^2 + \tau \varepsilon^{-2}C\dot{C}_n^2  \lVert \Delta_N \bm u^{n} \rVert_{2}^4.
  \end{aligned}
\end{equation}
Meanwhile, by applying the Sobolev interpolation inequality, we obtain
\[
  \left\Vert \Delta_N \bm u^n \right\Vert_2 \leq \left\Vert \nabla_N \bm u^n \right\Vert_2^{\frac{2}{3}}\left\Vert \Delta_N^2 \bm u^n \right\Vert_2^{\frac{1}{3}}.
\]
Consequently, applying Young's inequality to the term on the right-hand side of \eqref{eqn:stage2_h2_10} gives
\begin{equation}\label{eqn:stage2_h2_11}
\begin{aligned}
2(\kappa^2 + 3)\varepsilon^{-2} \left\Vert \Delta_N \bm u^n \right\Vert_2^2
\leq & 2(\kappa^2 + 3) \varepsilon^{-2}\left\Vert \nabla_N \bm u^n \right\Vert_2^{\frac{4}{3}}\left\Vert \Delta_N^2 \bm u^n \right\Vert_2^{\frac{2}{3}}  \\
\leq & C \dot{C}_n^2(1 + \kappa^3)\varepsilon^{-4}   + \frac{\varepsilon^{2}}{4}  \lVert \Delta_N^2 \bm u^n \rVert_{2}^2 \\
\end{aligned}
\end{equation}
and
\begin{equation}\label{eqn:stage2_h2_12}
\begin{aligned}
    \varepsilon^{-2} C \dot{C}_n^2 \lVert \Delta_N \bm u^n \rVert_{2}^4 \leq &  \varepsilon^{-2} C \dot{C}_n^2 \left\Vert \nabla_N \bm u^n \right\Vert_2^{\frac{8}{3}}\left\Vert \Delta_N^2 \bm u^n \right\Vert_2^{\frac{4}{3}} \\
    \leq &  C\dot{C}_n^{14}\varepsilon^{-10}  + \frac{\varepsilon^{2}}{4} \lVert \Delta_N^2 \bm u^n \rVert_{2}^2.
\end{aligned}
\end{equation}
Next, substituting \eqref{eqn:stage2_h2_11} and \eqref{eqn:stage2_h2_12} into \eqref{eqn:stage2_h2_10}, we obtain
\begin{equation}\label{eqn:stage2_h2_13}
    \begin{aligned}
    &\Vert (\phi_1 (\tau L_\kappa))^{-\frac{1}{2}} \Delta_N \bm u^{n+1} \Vert_2^2 + \frac{\varepsilon^2}{2} \tau \Vert \Delta_N^2 \bm u^{n+1} \Vert_2^2 + \frac{\kappa}{2} \left\Vert \nabla_N \Delta_N \bm u^{n+1} \right\Vert_2^2 \\
    \leq & \Vert (\phi_1 (\tau L_\kappa))^{-\frac{1}{2}} \Delta_N \bm u^n \Vert_2^2 - \frac{\varepsilon^2}{2} \tau \Vert \Delta_N^2 \bm u^n \Vert_2^2 - \frac{\kappa}{2} \left\Vert \nabla_N \Delta_N \bm u^{n} \right\Vert_2^2  +  \widehat{C}_{\varepsilon,1} \tau,
    \end{aligned}
\end{equation}
where we denote $\widehat{C}_{\varepsilon,1} =  C\dot{C}_n^2 (1 + \kappa^3)\varepsilon^{-4} + C\dot{C}_n^{14}\varepsilon^{-10}$. Moreover, by virtue of the Poincar\'e inequality stated in \eqref{eqn:poincareIneq}, for any $\bm f \in \mathcal{M}_N$, the following chain of inequalities holds
\begin{equation}
\begin{aligned}
     \left\lVert (\phi_1(\tau L_{\kappa}))^{-\frac{1}{2}} \Delta_N \bm f\right\rVert_{2}^2 \leq & \left\langle \frac{\tau L_{\kappa}}{1- e^{-\tau L_{\kappa}}} \Delta_N \bm f ,\Delta_N \bm f \right\rangle  \\
\leq & \left\langle (I + \tau L_{\kappa}) \Delta_N \bm f, \Delta_N \bm f \right\rangle \\
\leq & \lVert \Delta_N \bm f \rVert_{2}^2 + \tau (\varepsilon^2 \lVert \Delta_N^2 \bm f \rVert_{2}^2 + \kappa \lVert \nabla_N \Delta_N \bm f \rVert_{2}^2 ) \\
\leq & \frac{C + \varepsilon^2 \tau}{\varepsilon^2\tau} \tau (\varepsilon^2 \lVert \Delta_N^2 \bm f \rVert_{2}^2 + \kappa \lVert \nabla_N \Delta_N \bm f \rVert_{2}^2).
\end{aligned}
\label{Xeqn70-3.56}
\end{equation}
Returning to \eqref{eqn:stage2_h2_13}, we derive
\[
\begin{aligned}
         \frac{2C + 3\varepsilon^2 \tau}{2C + 2\varepsilon^2 \tau}\left\Vert (\phi_1 (\tau L_\kappa))^{-\frac{1}{2}} \Delta_N \bm u^{n+1} \right\Vert_2^2
    \leq&  \frac{2C + \varepsilon^2 \tau}{2C + 2\varepsilon^2 \tau} \left\Vert (\phi_1 (\tau L_\kappa))^{-\frac{1}{2}} \Delta_N \bm u^n \right\Vert_2^2 + \widehat{C}_{\varepsilon,1}\tau.
\end{aligned}
\]
Define $\widehat{C}_{\varepsilon,0} = C + \frac{1}{2}\varepsilon^2 \tau$. A combination of the recursive analysis with the a priori assumption reveals that
\begin{equation}\label{eqn:stage2_h2_est}
\resizebox{0.92\textwidth}{!}{$
\begin{aligned}
     \left\Vert (\phi_1 (\tau L_\kappa))^{-\frac{1}{2}} \Delta_N \bm u^{n+1} \right\Vert_2^2 \leq& \ddot{C}_{n+1}^2 := \Big(1- \frac{\varepsilon^2 \tau}{\widehat{C}_{\varepsilon,0} + \varepsilon^2 \tau}\Big)^{n+1} \ddot{C}_0^2 + C \Big(1- \Big(1- \frac{ \varepsilon^2 \tau}{\widehat{C}_{\varepsilon,0} + \varepsilon^2 \tau}\Big)^{n+1}\Big) (1 + \varepsilon^2\tau ) \varepsilon^{-2}\widehat{C}_{\varepsilon,1}.
\end{aligned}$}
\end{equation}
Here, $\ddot{C}_0$, $\widehat{C}_{\varepsilon,0}$ and $\widehat{C}_{\varepsilon,1}$ are global-in-time constants, which are independent of the final time $T$. Therefore, by invoking the a priori assumption \eqref{eqn:apriori_psi} at the previous time step, we are able to derive the $H^{2}$ norm estimate \eqref{eqn:stage2_h2_est}, which in turn validates the a priori assumption at the next time step.

\subsection{Recovery of the a priori assumption}\label{Xsec12-3.6}
The $\ell^\infty$ bound of numerical solutions $u_{n,1}$ has been derived in \eqref{eqn:stage1_H1refineEst} under the time step constraints
\[
    \tau \leq \tau_s := \min\{\tau_1, \tau_2\},
\]
where $\tau_1$ and $\tau_2$ are  defined in \eqref{eqn:stepConst_1} and \eqref{eqn:stepConst_2}, respectively. For the upper bound of $\Vert \bm u^{n+1} \Vert_\infty$, when $\varepsilon^2 \tau\leq 1$, an observation based on the uniform-in-time $H^2$ estimate yields that
\[
    \Vert \Delta_N u^{n+1}\Vert_2^2 \leq \Vert (\phi_1 (\tau L_\kappa))^{-\frac{1}{2}} \Delta_N \bm u^{n+1}\Vert_2^2 \leq \ddot{C}_{n+1}.
\]
In turn, by applying the Sobolev embedding inequality, we have
\[
    \Vert \bm u^{n+1}\Vert_{\infty} \leq C \Vert \nabla_N \bm u^{n+1} \Vert_2^{\frac{1}{2}} \Vert \Delta_N \bm u^{n+1} \Vert_2^{\frac{1}{2}} \leq \sqrt{2} \tilde{C}_{n}.
\]
Furthermore, if we select the fixed constant $\kappa$ as specified in \linkref[Theorem]{\ref{thm:p_stab}}, and impose the time-step constraint $\tau \leq \tau_s :=  \min\{\tau_1,\tau_2\}$, the energy stability of the next step is thereby established, i.e.,
\[
    E_N(\bm u^{n+1}) \leq E_N(\bm u^n) \leq E_N(\bm u^0) := C_e.
\]
This verifies the a priori assumption at the next time step. Consequently, the theoretical result of this work has been successfully proven. Meanwhile, with the aid of \linkref[Lemma]{\ref{lem:energy_est}}, the uniform-in-time bounds for $H^1$ and $\ell^\infty$ at the next time step are also justified, which are given by
\[
    \Vert \nabla_N \bm u^{n+1}\Vert_2 \leq \dot{C}_n := \varepsilon^2 (2C_e)^{\frac{1}{2}},\quad \Vert \bm u^{n+1} \Vert_\infty \leq \tilde{C}_n.
\]
Thus, the induction argument presented in this section can be effectively applied.

\begin{remark}    To streamline the presentation of both the numerical method and the theoretical analysis, periodic boundary conditions are adopted in this work. This choice notably enhances the efficiency of implementing the ERK2 scheme through the Fourier spectral collocation method and simplifies the subsequent analysis.
    For boundary conditions with more physical relevance, such as homogeneous Neumann or Dirichlet boundary conditions, the proposed scheme and analytical framework can be adapted without significant difficulty. In particular, under homogeneous Neumann boundary conditions, the phase variable $ u$ and its Laplacian $\Delta u$ can be approximated using the Fourier cosine transform, which preserves spectral accuracy. As a result, all Fourier-based operators can be defined in a fully analogous manner, and the corresponding stability and energy estimates can be derived following the same line of reasoning. This enables the establishment of global-in-time energy stability for the ERK2 scheme under these boundary conditions as well.
\label{Xenun5-2}
\end{remark}

\section{Convergence analysis}\label{Xsec13-4}\label{sec:error}
By using the uniform boundedness of the numerical solution, as well as the preliminary operator estimates stated in \linkref[Section]{\ref{sec:stab}}, we are able to derive the $\ell^2$ error estimate for the ERK2 numerical scheme to conclude this work.

Let the exact solution of the CH \linkref[\del{equation}\ins{Eq.}]{\eqref{eqn:ch_eq}} be denoted by $u_e$, and assume that it satisfies
\[
    u_e\in \mathcal{R} :=  L^2(0,T; H_{\text{per}}^{m_0 + 2}(\Omega)) \cap L^{\infty} (0,T; H_\text{per}^{m_0}(\Omega)),
\]
with sufficiently smooth initial data. Define $U_N(x,y,z,t):= \mathcal{P}_N u_e(x,y,z,t)$, where $\mathcal{P}_N$ denotes the Fourier projection operator onto $\mathcal{M}_N$. The following projection estimate is standard: if $u_e \in L^\infty(0,T; H^m_{per}(\mathrm{\Omega}))$ for $m \geq k \geq 0$, then
\[
    \Vert U_N - u_e \Vert_{L^{\infty}(0,T: H^k)} \leq Ch^{m-k}  \Vert u_e \Vert_{L^{\infty}(0,T: H^m)}.
\]

In addition, let $U^n$ be the grid values of the projection solution $U_N$ restricted on the mesh $\Omega_N$ at fixed time $t_n$. Then $U^n \in \mathcal{M}_N$, and for all $(p,q,r)\in \mathcal{S}_N$, it holds that $U^n_{p,q,r} = U_N(x_p,y_q,z_r,t_n)$. The initial data are specified by $U_{p,q,r}^0 = U_N(x_p,y_q,z_r,0)$.

We first introduce a reference function $U_{n,1}$. By performing a careful consistency analysis with the aid of \linkref[Lemma]{\ref{lem:erk2_err}}, we derive that the following relations hold:
\begin{align}
    U_{n,1} &= \phi_0(\tau L_\kappa) U^n + \tau \phi_1(\tau L_\kappa) N_\kappa \left(U^n\right), \label{eqn:consisStage_1}\\
    U^{n+1} &= \phi_0(\tau L_\kappa) U^n + \frac{\tau}{2} \phi_1(\tau L_\kappa) (N_\kappa (U^n) + N_\kappa (U_{n,1})) + \tau \zeta^n, \label{eqn:consisStage_2}
\end{align}
with $\Vert\lap_N \zeta\Vert_2 = \mathcal{O}(\tau^2 + h^{m_0})$. Next, we define the error grid functions as follows:
\[
    \bm e^{n} = U^n - \bm u^n,\quad \bm e_{n,1} = U_{n,1} - \bm u_{n,1}, \quad \forall n \geq 0.
\]
Subtracting the numerical schemes \eqref{eqn:etd_s1}--\eqref{eqn:etd_s2} from  \eqref{eqn:consisStage_1}--\eqref{eqn:consisStage_2} yields
\begin{align}
    &\bm e_{n,1} = \phi_0(\tau L_\kappa) \bm e^n + \tau \phi_1(\tau L_\kappa) (N_\kappa \left(U^n\right) - N_\kappa (\bm u^n) ),\label{eqn:errEq_stage1}\\
    &\bm e^{n+1} = \phi_0(\tau L_\kappa) \bm e^n + \frac{\tau}{2} \phi_1(\tau L_\kappa) (N_\kappa \left(U^n\right) - N_\kappa (\bm u^n) +N_\kappa (U_{n,1}) - N_\kappa (\bm u_{n,1}) ) + \tau \zeta^n.\label{eqn:errEq_stage2}
\end{align}
We now state the convergence estimate for the proposed fully discrete scheme.
\begin{theorem}\label{thm:conv}
 Given initial data $u^0 \in C_{\text{per}}^{m_0+2}(\Omega)$ with periodic boundary conditions, suppose that the unique solution to the CH equation belongs to the regularity class $\mathcal{R}$. For the ERK2 scheme, if $\tau$ and $h$ are sufficiently small, the solution of the fully discrete scheme satisfies the following convergence estimate:
 \[
     \Vert \bm e^n\Vert_2 \leq C(\tau^2 + h^{m_0}),\quad \forall n \leq N_t,
 \]
 where $C>0$ is a constant independent of $\tau$ and $h$.
\end{theorem}

\begin{proof}
Following an argument analogous to the stability analysis, the error \linkref[\del{equations}\ins{Eqs.}]{\eqref{eqn:errEq_stage1}}--\eqref{eqn:errEq_stage2} can be reformulated as
\begin{align}
&\bm e_{n,1} - \bm e^{n} + \tau\phi_1(\tau L_\kappa) \bm e^n = \tau \phi_1(\tau L_\kappa) \left(N_\kappa \left(U^n\right) - N_\kappa (\bm u^n) \right),\label{eqn:errEq_stage1r}\\[.5em]
&\bm e^{n+1} - \bm e^{n} + \tau\phi_1(\tau L_\kappa) \bm e^n = \frac{\tau}{2} \phi_1(\tau L_\kappa) \left(N_\kappa \left(U^n\right) - N_\kappa (\bm u^n) + N_\kappa (U_{n,1}) - N_\kappa (\bm u_{n,1}) \right) + \tau \zeta^n.\label{eqn:errEq_stage2r}
\end{align}
For the first stage, taking the discrete $\ell^2$-inner product of both sides of \eqref{eqn:errEq_stage1r}  with $\bm e_{n,1}$ yields
\begin{equation}\label{eqn:err_est_01}
\begin{aligned}
    &\Vert\bm e_{n,1} \Vert_2^2 - \Vert \bm e^{n} \Vert_2^2 + \Vert \bm e_{n,1} - \bm e^{n} \Vert_2^2 + \tau (\Vert G_N^{(0)} \bm e_{n,1}  \Vert^2 + \Vert G_N^{(0)} \bm e^{n} \Vert^2 - \Vert G_N^{(0)}(\bm e_{n,1} - \bm e^{n}) \Vert_2^2) \\
    = & 2\tau \langle  \phi_1(\tau L_\kappa)(N_\kappa \left(U^n\right) - N_\kappa (\bm u^n)), \bm e_{n,1} \rangle,
\end{aligned}
\end{equation}
where $G_N^{(0)} := (\phi_1(\tau L_\kappa) L_\kappa)^{\frac{1}{2}}$.
Furthermore, analogous to the operators introduced in \eqref{eqn:linear_operators}, $G_N^{(0)}$ inherits the properties stated in \linkref[Lemma]{\ref{lem:operator_est}}. Consequently, we have
\begin{equation}\label{eqn:err_01}
    \Vert \bm e_{n,1} - \bm e^{n} \Vert_2^2 - \tau \Vert G_N^{(0)}(\bm e_{n,1} - \bm e^{n}) \Vert_2^2 \geq 0
\end{equation}
and
\begin{equation}\label{eqn:err_02}
\begin{aligned}
    \Vert G_N^{(0)} \bm e_{n,1} \Vert_2^2 + \Vert G_N^{(0)} \bm e^{n} \Vert_2^2
    =& \varepsilon^2 (\Vert \left( \phi_1(\tau L_\kappa)\right)^{\frac{1}{2}} \Delta_N \bm e_{n,1} \Vert_2^2 + \Vert \left( \phi_1(\tau L_\kappa)\right)^{\frac{1}{2}} \Delta_N \bm e^{n} \Vert_2^2) \\ &+ \kappa (\Vert \left( \phi_1(\tau L_\kappa)\right)^{\frac{1}{2}} \nabla_N \bm e_{n,1} \Vert_2^2  + \Vert \left( \phi_1(\tau L_\kappa)\right)^{\frac{1}{2}} \nabla_N \bm e^{n} \Vert_2^2).
\end{aligned}
\end{equation}
On the other hand, the inner product on the right-hand side of \eqref{eqn:err_est_01} can be written as
\[
    2\langle \phi_1(\tau L_\kappa)(N_\kappa \left(U^n\right) - N_\kappa (\bm u^n)), \bm e_{n,1} \rangle = 2\langle \phi_1(\tau L_\kappa)\Delta_N(\left(U^n\right)^3  -  (\bm u^n)^3), \bm e_{n,1} \rangle  - 2(\kappa+1)\langle \phi_1(\tau L_\kappa)\Delta_N\bm e^n, \bm e_{n,1} \rangle.
\]
By applying the estimate in \eqref{eqn:phi_esitmate} and introducing the notation $\tilde{C} = \max\{\Vert U^n\Vert_\infty, \Vert  \bm u^n \Vert_\infty \}$, we obtain
\begin{equation}\label{eqn:err_03}
    2\langle \phi_1(\tau L_\kappa)\Delta_N(\left(U^n\right)^3 -  (\bm u^n)^3), \bm e_{n,1} \rangle \leq \frac{6\tilde{C}^2}{\varepsilon^2} \Vert \bm e^n \Vert^2  + \frac{\varepsilon^2}{2}\Vert \left( \phi_1(\tau L_\kappa)\right)^{\frac{1}{2}}\Delta_N \bm e_{n,1} \Vert_2^2
\end{equation}
and
\begin{equation}\label{eqn:err_04}
\begin{aligned}
     -2(\kappa+1)\langle \phi_1(\tau L_\kappa)\Delta_N\bm e^n, \bm e_{n,1} \rangle \leq &\kappa(\Vert \left( \phi_1(\tau L_\kappa)\right)^{\frac{1}{2}} \nabla_N \bm e^n \Vert_2^2 + \Vert \left( \phi_1(\tau L_\kappa)\right)^{\frac{1}{2}} \nabla_N \bm e_{n,1} \Vert_2^2) \\
     &+ \frac{2}{\varepsilon^2} \Vert \bm e^{n} \Vert_2^2 + \frac{\varepsilon^2}{2} \Vert \left( \phi_1(\tau L_\kappa)\right)^{\frac{1}{2}}\Delta_N \bm e_{n,1} \Vert_2^2.
\end{aligned}
\end{equation}
Substituting \eqref{eqn:err_01}--\eqref{eqn:err_04} into \eqref{eqn:err_est_01} gives
\begin{equation}
    \Vert\bm e_{n,1} \Vert_2^2 - \Vert \bm e^{n} \Vert_2^2 + \tau \varepsilon^2 \Vert \left( \phi_1(\tau L_\kappa)\right)^{\frac{1}{2}}\Delta_N \bm e^{n} \Vert_2^2 \leq \tau \frac{2 + 6\tilde{C}^2}{\varepsilon^2}\Vert \bm e^{n} \Vert_2^2.
\label{Xeqn77-4.12}
\end{equation}
As an immediate consequence of the above inequality, we establish a preliminary $\ell^2$ error estimate for the first stage of the ERK2 scheme:
\begin{equation}\label{eqn:err_05}
    \Vert\bm e_{n,1} \Vert_2^2 \leq 2 \Vert \bm e^{n} \Vert_2^2,
\end{equation}
provided that the time step satisfies the constraint $\tau \leq \varepsilon^2/(1+3\tilde{C}^2)$.

For the second \linkref[\del{equation}\ins{Eq.}]{\eqref{eqn:errEq_stage2r}}, taking the discrete $\ell^2$ inner product with $2 \bm e^{n+1}$ and using the inequality \eqref{eqn:err_01} yields
\begin{equation}\label{eqn:err_est_02}
\begin{aligned}
   &\Vert\bm e^{n+1} \Vert_2^2 - \Vert \bm e^{n} \Vert_2^2 + \tau (\Vert G_N^{(0)} \bm e^{n+1}  \Vert^2 + \Vert G_N^{(0)} \bm e^{n} \Vert^2) \\
   \leq & \tau \langle \phi_1(\tau L_\kappa)(N_\kappa \left(U^n\right) - N_\kappa (\bm u^n)), \bm e^{n+1} \rangle + \tau \langle  \phi_1(\tau L_\kappa)(N_\kappa (U_{n,1}) - N_\kappa (\bm u_{n,1})), \bm e^{n+1} \rangle + 2\tau \langle \zeta^n , \bm e^{n+1}\rangle.
\end{aligned}
\end{equation}
The first nonlinear inner product can be analyzed in the same way as in \eqref{eqn:err_03}--\eqref{eqn:err_04}, and thus, we obtain
\begin{equation}\label{eqn:err_06}
\begin{aligned}
    &\langle \phi_1(\tau L_\kappa)(N_\kappa \left(U^n\right) - N_\kappa (\bm u^n)), \bm e^{n+1} \rangle \\
    \leq & \frac{3\tilde{C}^2}{\varepsilon^2} \Vert \bm e^n \Vert^2  + \frac{\varepsilon^2}{4}\Vert \left( \phi_1(\tau L_\kappa)\right)^{\frac{1}{2}}\Delta_N \bm e^{n+1} \Vert_2^2 + \frac{\kappa}{2}(\Vert \left( \phi_1(\tau L_\kappa)\right)^{\frac{1}{2}} \nabla_N \bm e^n \Vert_2^2 + \Vert \left( \phi_1(\tau L_\kappa)\right)^{\frac{1}{2}} \nabla_N \bm e^{n+1} \Vert_2^2) \\
     &+ \frac{1}{\varepsilon^2} \Vert \bm e^{n} \Vert_2^2 + \frac{\varepsilon^2}{4} \Vert \left( \phi_1(\tau L_\kappa)\right)^{\frac{1}{2}}\Delta_N \bm e^{n+1} \Vert_2^2.
\end{aligned}
\end{equation}
For the second nonlinear inner product in \eqref{eqn:err_est_02}, we combine it with the estimate \eqref{eqn:err_05} derived from the first stage, which leads to
\begin{equation}\label{eqn:err_07}
    \begin{aligned}
        &\langle \phi_1(\tau L_\kappa)(N_\kappa (U_{n,1}) - N_\kappa (\bm u_{n,1})), \bm e^{n+1} \rangle \\
    \leq & \frac{3\tilde{C}^2}{\varepsilon^2} \Vert \bm e_{n,1} \Vert^2  + \frac{\varepsilon^2}{4}\Vert \left( \phi_1(\tau L_\kappa)\right)^{\frac{1}{2}}\Delta_N \bm e^{n+1} \Vert_2^2 + \frac{(\kappa+1)^2}{\varepsilon^2} \Vert \bm e_{n,1} \Vert_2^2 + \frac{\varepsilon^2}{4} \Vert \left( \phi_1(\tau L_\kappa)\right)^{\frac{1}{2}}\Delta_N \bm e^{n+1} \Vert_2^2\\
    \leq & \frac{6\tilde{C}^2}{\varepsilon^2} \Vert \bm e^{n} \Vert^2  + \frac{\varepsilon^2}{4}\Vert \left( \phi_1(\tau L_\kappa)\right)^{\frac{1}{2}}\Delta_N \bm e^{n+1} \Vert_2^2 + \frac{2(\kappa+1)^2}{\varepsilon^2} \Vert \bm e^{n} \Vert_2^2 + \frac{\varepsilon^2}{4} \Vert \left( \phi_1(\tau L_\kappa)\right)^{\frac{1}{2}}\Delta_N \bm e^{n+1} \Vert_2^2.
    \end{aligned}
\end{equation}
Furthermore, the bound for the truncation error inner product term can be derived straightforwardly as follows:
\begin{equation}\label{eqn:err_08}
    2 \langle \zeta^n , \bm e^{n+1}\rangle \leq \Vert \zeta^n \Vert_2^2 + \Vert \bm e^{n+1} \Vert_2^2.
\end{equation}
Consequently, by combining \eqref{eqn:err_est_02}--\eqref{eqn:err_08}, we obtain
\[
\begin{aligned}
    &\Vert\bm e^{n+1} \Vert_2^2 - \Vert \bm e^{n} \Vert_2^2  + \tau\varepsilon^2 \Vert \left( \phi_1(\tau L_\kappa)\right)^{\frac{1}{2}}\Delta_N \bm e^{n} \Vert_2^2 + \frac{\tau\kappa}{2}(\Vert \left( \phi_1(\tau L_\kappa)\right)^{\frac{1}{2}} \nabla_N \bm e_{n,1} \Vert_2^2  + \Vert \left( \phi_1(\tau L_\kappa)\right)^{\frac{1}{2}} \nabla_N \bm e^{n} \Vert_2^2) \\
    \leq &\tau \frac{9\tilde{C}^2 +1 + 2 (\kappa+1)^2}{\varepsilon^2} \Vert \bm e^{n}  \Vert^2 +  \tau\Vert \bm e^{n+1} \Vert_2^2 + \tau\Vert \zeta^n \Vert_2^2.
\end{aligned}
\]
Summing the above inequality from $0$ to $n$ with $\bm e^0 = 0$ gives
\[
\begin{aligned}
    \Vert\bm e^{n+1} \Vert_2^2 + \tau \varepsilon^2 \sum_{k = 0}^{n} \Vert \left( \phi_1(\tau L_\kappa)\right)^{\frac{1}{2}}\Delta_N \bm e^{k} \Vert_2^2 \leq & \tau  \frac{9\tilde{C}^2 + 1 + \varepsilon^2  + 2 (\kappa+1)^2}{\varepsilon^2} \sum_{k=0}^{n} \Vert \bm e^{k} \Vert_2^2   + \tau \Vert \bm e^{n+1} \Vert_2^2 + T \sup_{k=0,\dots,n} \Vert \zeta^{k} \Vert_2^2. \\
    \leq & \tau  \frac{9\tilde{C}^2 + 1 + \varepsilon^2  + 2 (\kappa+1)^2}{\varepsilon^2} \sum_{k=0}^{n} \Vert \bm e^{k} \Vert_2^2   + \tau \Vert \bm e^{n+1} \Vert_2^2 + T (\tau^2 + h^{m_0})^2.
\end{aligned}
\]
Let
\[C_k = (9\tilde{C}^2 + 1 + \varepsilon^2  + 2 (\kappa+1)^2)/\varepsilon^2,\]
and assume the time step satisfies
\[\tau \leq \min\{\frac{1}{2},\varepsilon^2/(1+3\tilde{C}^2)\}.\]
Under this condition, an application of the discrete Gr\"onwall inequality leads to
\[
    \Vert\bm e^{n+1} \Vert_2^2 + \tau \varepsilon^2 \sum_{k = 0}^{n} \Vert \left( \phi_1(\tau L_\kappa)\right)^{\frac{1}{2}}\Delta_N \bm e^{k} \Vert_2^2 \leq 2T\exp(2C_kT)(\tau^2 + h^{m_0})^2.
\]
This completes the proof.
\end{proof}

{\opecolor{blue}
\section{Numerical experiments}\label{Xsec14-5}\label{sec:numex}
In this section, several numerical experiments are conducted to verify the theoretical results of the proposed scheme \eqref{eqn:etd_s1}--\eqref{eqn:etd_s2}. We will verify the temporal convergence rate and demonstrate the robustness of the scheme over a long-time simulation of phase separation in three-dimensional space. Theoretically, the stabilization constant $\kappa$ and time step size $\tau$ need to satisfy the condition presented in \linkref[Theorem]{\ref{thm:p_stab}} to ensure energy stability. In practice, the following experiments will demonstrate that $\kappa = 2$ is already sufficient. In all parts below, the spatial discretization is always performed using the Fourier pseudo-spectral method with the $2/3$-rule for dealiasing on a periodic cubic domain.
{\opecolor{blue}
\begin{example}\label{ex:stab_test}
Consider
\begin{equation}
\nonumber
    \left\{
    \begin{aligned}
    &u_t = \Delta (- \varepsilon^2 \Delta u +  f(u)),\quad \bm{x} \in (0,2\pi)^3,\ t \in (0,20],\\
    &u(\bm{x},0) = u_0(\bm{x}),
    \end{aligned}
    \right.
\end{equation}
with $\varepsilon=0.01$ and $ u_0(\bm{x})$ generated by a zero-mean random field with amplitude $0.5$.
\end{example}

In this example, we investigate the long-time stability behavior of the stabilized ERK2 scheme \eqref{eqn:etd_s1}--\eqref{eqn:etd_s2} with 128 Fourier modes used in space for the 3D CH \linkref[\del{equation}\ins{Eq.}]{\eqref{eqn:ch_eq}}. We first compute a reference solution with the small time step size $\tau=10^{-3}$ and  $\kappa=0$. \linkref[\del{Figure}\ins{Fig.}]{\ref{fig:ch3d_small_tau_kappa0_norms}} reports the temporal evolution of the discrete $H^1$ semi-norm, the discrete $H^2$ semi-norm, and the $\ell^\infty$ norm of the numerical solution. The result in \linkref[\del{Figure}\ins{Fig.}]{\ref{fig:ch3d_small_tau_kappa0_norms}}   shows that the numerical solution remains uniformly bounded throughout the simulation. This observation is consistent with the uniform-in-time $H^1$, $H^2$, and $\ell^\infty$ bounds established in \linkref[Theorem]{\ref{thm:long_time_stability}}, and it provides a reliable baseline for assessing the behavior of large-step computations.
\begin{figure}[t]
\begin{auimage}{}
\centering
\includegraphics[width=0.9\linewidth]{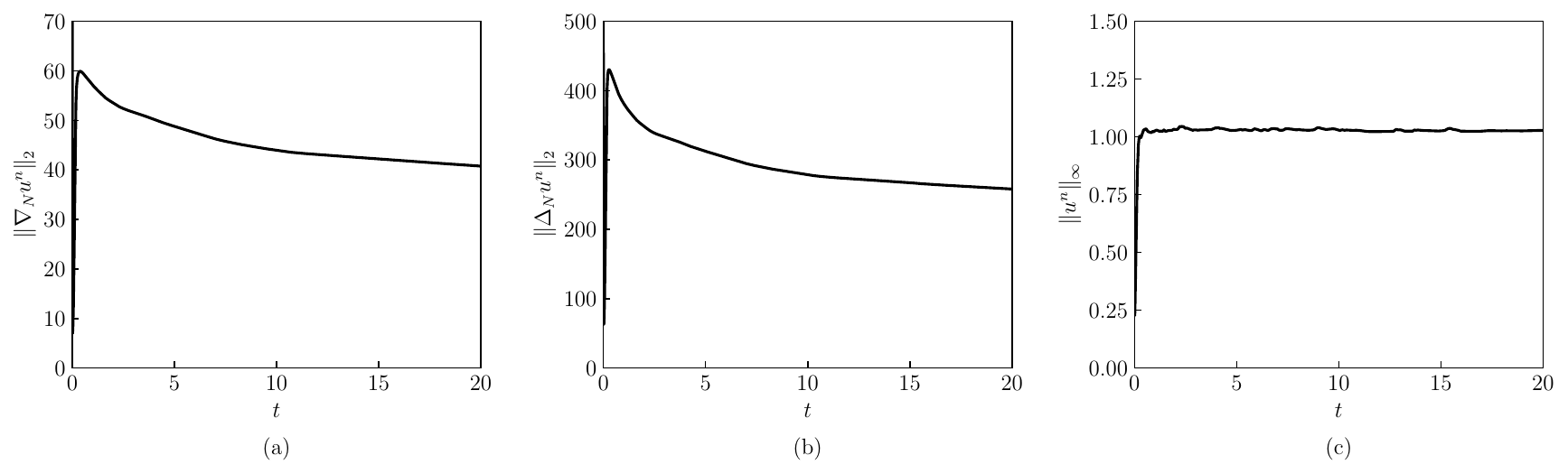}
\end{auimage}

\caption{Temporal evolution of the discrete $H^1$ semi-norm, discrete $H^2$ semi-norm, and $\ell^\infty$ norm for the reference solution of the ERK2 scheme \eqref{eqn:etd_s1}--\eqref{eqn:etd_s2} with $\tau=10^{-3}$ and $\kappa=0$ for \linkref[Example]{\ref{ex:stab_test}}.\label{fig:ch3d_small_tau_kappa0_norms}}
\end{figure}
\begin{figure}[t]
\begin{auimage}{}
\centering
\includegraphics[width=0.9\linewidth]{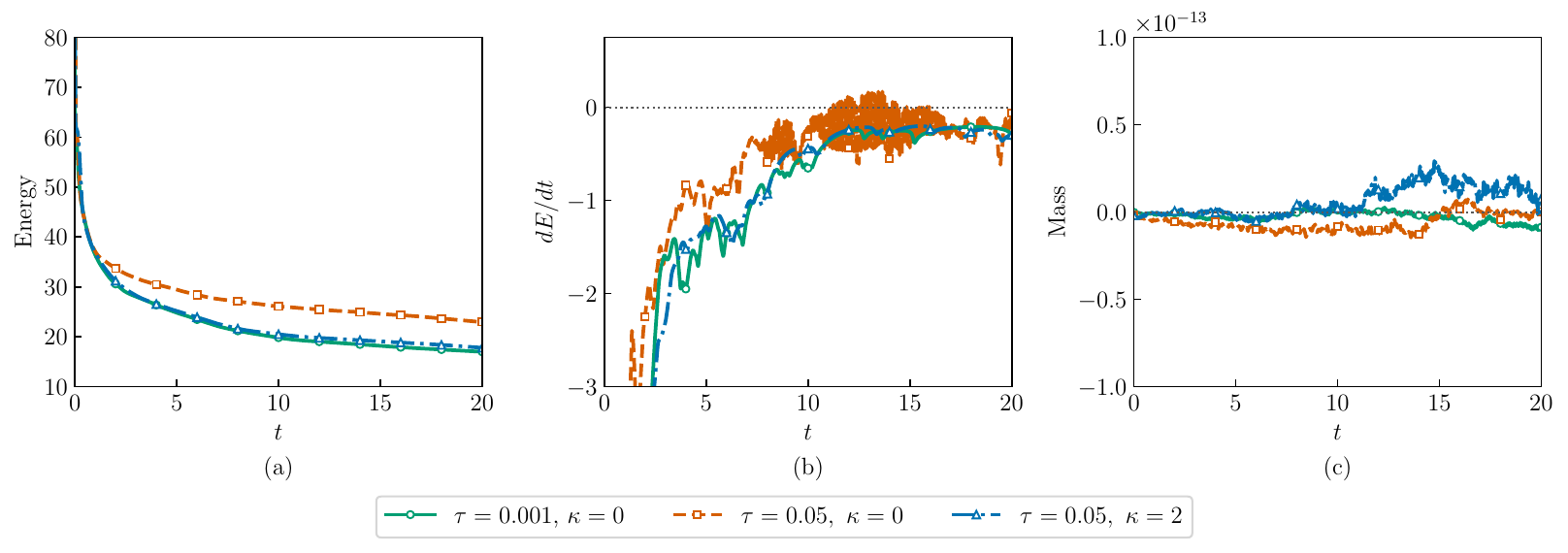}
\end{auimage}

\caption{Energy, energy dissipation rate, and mass profiles computed by the ERK2 scheme \eqref{eqn:etd_s1}--\eqref{eqn:etd_s2} for \linkref[Example]{\ref{ex:stab_test}}. The reference solution uses $\tau=10^{-3}$ and $\kappa=0$, while the large-step computations use $\tau=5\times 10^{-2}$ for both the unstabilized case ($\kappa=0$) and the stabilized case ($\kappa=2$).\label{fig:ch3d_big_tau_energy_contrast}}
\end{figure}

Subsequently, we compare the reference solution with large time step computations using different stabilization choices. As shown in \linkref[\del{Figure}\ins{Fig.}]{\ref{fig:ch3d_big_tau_energy_contrast}}, the reference solution corresponding to $\tau=10^{-3}$ and $\kappa=0$ is plotted alongside two numerical results obtained with a larger time step $\tau=5\times 10^{-2}$: the unstabilized scheme with $\kappa=0$ and the stabilized scheme with $\kappa=2$. The time derivative curve of energy in the profile indicates that increasing the time step size may cause the method to deviate from the anticipated dissipative characteristics. In contrast, stabilized computation eliminates such numerical instability and recovers  energy dissipation behavior. This comparison demonstrates the practical effect of the stabilization parameter within the energy stability criterion presented in \linkref[Theorem]{\ref{thm:p_stab}}. Meanwhile, the mass distribution stays nearly at the machine rounding error magnitude across all numerical cases, verifying that the Fourier pseudo-spectral discretization and time marching scheme preserve the conservative properties of CH dynamics within numerical precision bounds.
\begin{figure}[t]
\begin{auimage}{}
\centering
\includegraphics[width=0.65\linewidth]{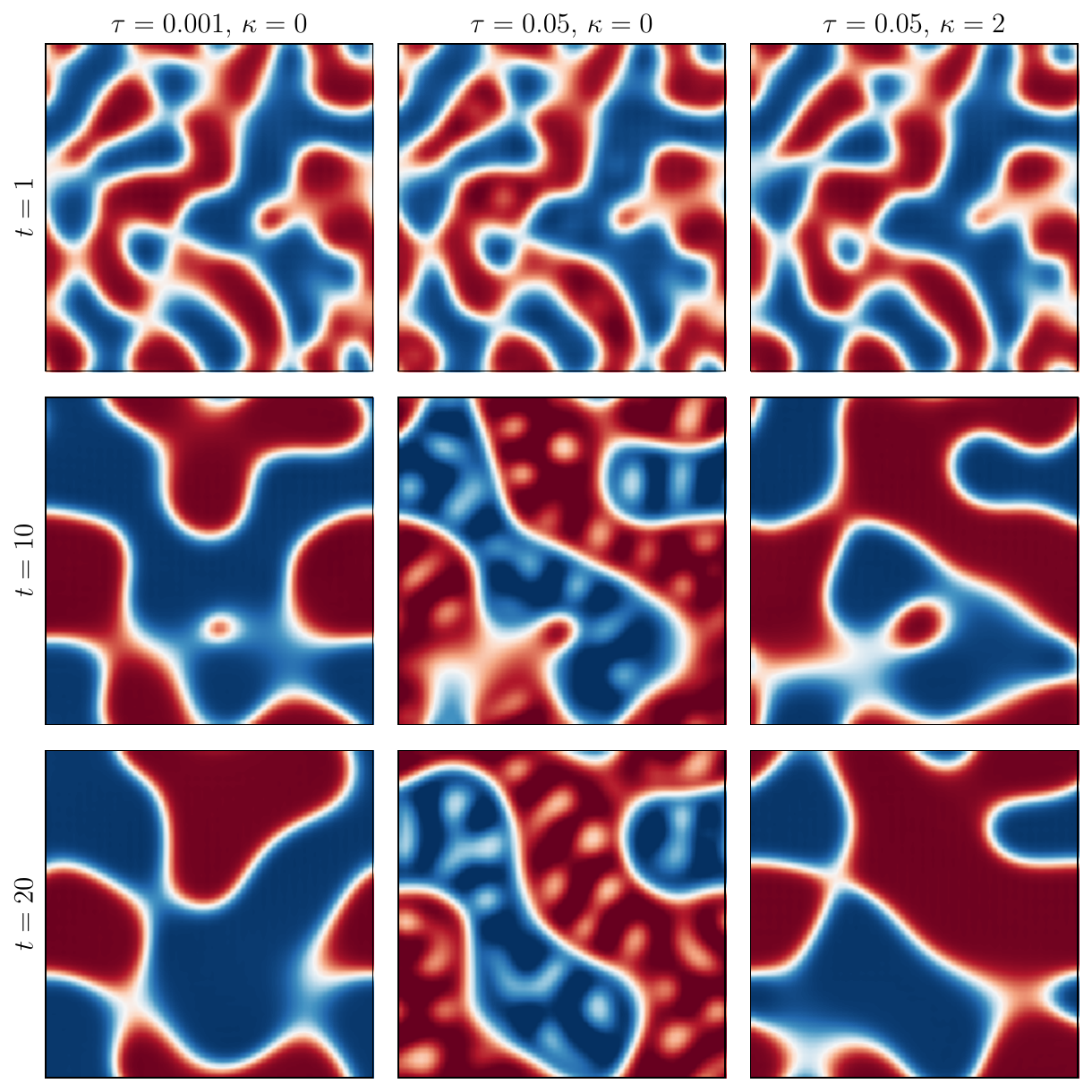}
\end{auimage}

\caption{Slices $u(x,y,\pi,t)$ of the three-dimensional phase field at $t=1$, $10$, and $20$ for \linkref[Example]{\ref{ex:stab_test}}, computed by the ERK2 scheme \eqref{eqn:etd_s1}--\eqref{eqn:etd_s2} using different choices of the time step size and stabilization parameter.\label{fig:ch3d_phase_snapshots}}
\end{figure}
}

Finally, \linkref[\del{Figure}\ins{Fig.}]{\ref{fig:ch3d_phase_snapshots}} displays two-dimensional slices of the three-dimensional phase field at $z=\pi$ for $t=1$, 10 and 20. The first column corresponds to the reference computation, while the remaining columns show the large-step unstabilized and stabilized computations. The unstabilized large-step solution in the second column shows visibly distorted phase patterns compared to the reference solution, consistent with the loss of energy dissipation observed in \linkref[\del{Figure}\ins{Fig.}]{\ref{fig:ch3d_big_tau_energy_contrast}}. By contrast, the stabilized solution in the third column eliminates these false oscillatory characteristics and yields phase profiles whose qualitative evolution better matches the reference behavior. These results numerically support the main theoretical conclusions in \linkref[Theorems]{\ref{thm:p_stab}} and~\ref{thm:long_time_stability}: the ERK2 scheme  \eqref{eqn:etd_s1}--\eqref{eqn:etd_s2} preserves mass to numerical precision, exhibits uniform-in-time boundedness in the monitored norms for sufficiently small time steps, and benefits from stabilization in maintaining the dissipative property in simulations.

\begin{example}\label{exm:ch3d_conv_test}
In this example, we test the temporal convergence of the proposed second-order scheme \eqref{eqn:etd_s1}--\eqref{eqn:etd_s2} using the method of manufactured solutions.
Set $\Omega = (0,2\pi)^3$, $\varepsilon^2 = 0.01$, and $\kappa = 2$.
A smooth exact solution is prescribed as
\begin{equation}\label{eqn:mms_ch3d}
u_e(x,y,z,t) = 0.4\cos(t)\cos(x)\cos(y)\cos(z).
\end{equation}
The CH \linkref[\del{equation}\ins{Eq.}]{\eqref{eqn:ch_eq}} is augmented with an external forcing term $g$, which is determined by substituting \eqref{eqn:mms_ch3d}.
This external forcing term is also incorporated into the numerical scheme \eqref{eqn:etd_s1}--\eqref{eqn:etd_s2} when computing the numerical solution, allowing the numerical error to be accurately evaluated.
\end{example}

For the numerical results shown in \linkref[\del{Figure}\ins{Fig.}]{\ref{fig:convergency_ch3d}}, uniform time steps $\tau = 0.01 \times 2^{-k}$ with $k = 2, 3, \dots, 6$ are employed in the temporal discretization.
To ensure that spatial truncation errors remain negligible relative to temporal errors, $64$ Fourier modes are used in each spatial direction.
At the final time $T = 1$, the discrete $\ell^{\infty}$, $\ell^{2}$ and $H^{1}$ errors between the numerical solution and $u_e$ from \eqref{eqn:mms_ch3d} are computed. As can be seen from \linkref[\del{Figure}\ins{Fig.}]{\ref{fig:convergency_ch3d}}, the proposed ERK2 scheme \eqref{eqn:etd_s1}--\eqref{eqn:etd_s2} remains stable for all tested time steps. The empirical convergence rates for the $\ell^{\infty}$, $\ell^{2}$, and $H^{1}$ errors all approach $2$, confirming the expected second-order temporal convergence of the scheme in all three norms. These results are consistent with the theoretical temporal accuracy established in \linkref[Theorem]{\ref{thm:conv}}.
\begin{figure}[t]
\begin{auimage}{}
\centering
\includegraphics[width=0.6\linewidth]{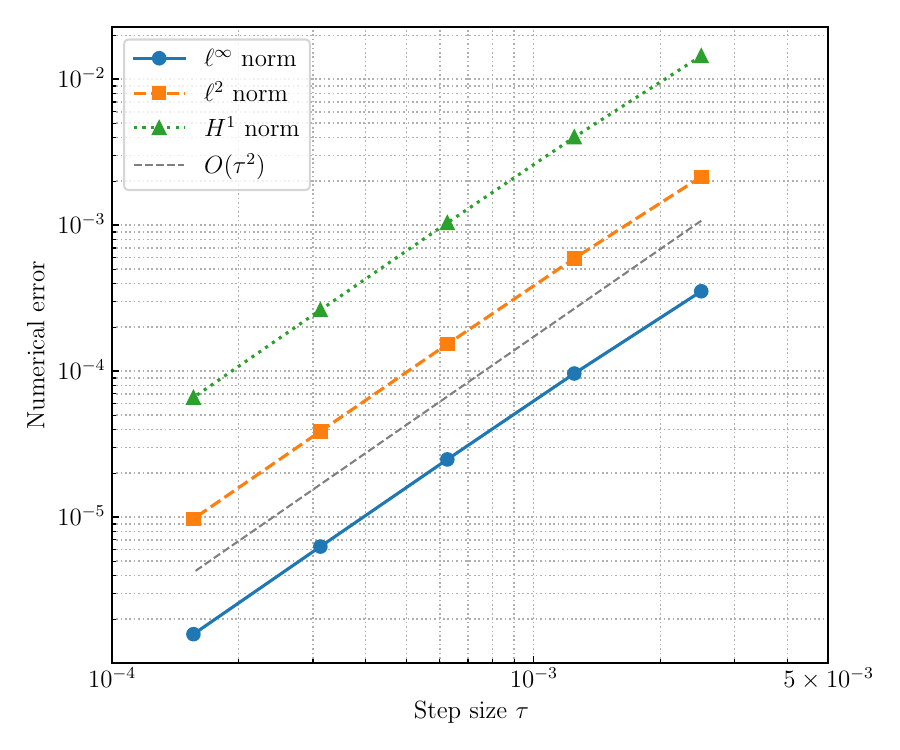}
\end{auimage}

\caption{The discrete $\ell^{\infty}$, $\ell^{2}$, and $H^{1}$ errors at $T = 1$ between the numerical solution computed by the proposed second-order ERK2 scheme and the manufactured exact solution $u_e$ in
\eqref{eqn:mms_ch3d}, plotted against the time-step size $\tau = 0.01 \times 2^{-k}$, for $k = 2, 3, \dots, 6$. The dashed gray line indicates the reference slope $\mathcal{O}(\tau^{2})$.\label{fig:convergency_ch3d}}
\end{figure}}

\section{Conclusion}\label{Xsec15-6}\label{sec:con}
This work provides a detailed analysis of a second-order accurate exponential Runge--Kutta (ERK2) numerical scheme for the Cahn--Hilliard equation. Utilizing the Sobolev embedding inequality, certain preliminary operator estimates, and an a priori assumption on the numerical solution from the previous time step and the initial data, we have rigorously established the long-time numerical stability of the scheme in both the $H^1$ and $H^2$ norms. Consequently, energy stability emerges as a direct corollary of this analysis. The proven uniform boundedness of the numerical solution further facilitates the derivation of an optimal error estimate for the ERK2 scheme. {The theoretical results are validated by numerical experiments, including  long-time three-dimensional simulations and a convergence test.}  The methodological framework developed here for long-time stability analysis can be extended to a broader class of high-order, multi-stage Runge--Kutta-type exponential integrators and is applicable to a wide range of related problems.

\begin{coi}
\section*{Declaration of competing interest}
 The authors declare that they have no known competing financial interests or personal relationships that could have appeared to influence the work reported in this paper.
\end{coi}

\begin{ack}
\section*{Acknowledgment}
The authors would like to thank the anonymous referees for the valuable comments and constructive suggestions that have led to significant improvements in this work.
\end{ack}

\section*{Data availability}
No data was used for the research described in the article.

\begin{appgroup}
\appsection{Proof of \linkref[Lemma]{\ref{lem:perFun_prop}}}\label{apps1}
Denote
$\bm{k} = (k, \ell, m) \in \widehat{S}_N$. 
 For the inequality \eqref{eqn:3dSob}, let $R$ be an arbitrary positive constant. By using the Fourier expansion of $\bm f$, one has
\begin{equation}\label{eqn:infexpend}
\frac{1}{3} \left\Vert \bm f \right\Vert_\infty^{2} \leq
\frac{1}{3} \left( \sum_{\bm k \in \widehat{\mathcal{S}}_N} \left\vert \hat{f}_{\bm k} \right\vert \right)^{\!2} \leq
\left\vert \hat{f}_{0,0,0} \right\vert^2 +
\left( \sum_{0 < \left\vert \bm k \right\vert \leq R} \left\vert \hat{f}_{\bm k} \right\vert \right)^{\!2} +
\left( \sum_{\left\vert \bm k \right\vert > R} \left\vert \hat{f}_{\bm k} \right\vert \right)^{\!2}.
\end{equation}
By the Cauchy--Schwarz inequality, the second and third terms on the right-hand side of \eqref{eqn:infexpend} can be estimated as
\[
\begin{aligned}
\left(\sum_{0 < |\bm k| \le R} \left|\hat{f}_{\bm {k}}\right| \right)^2 &\leq \left(\sum_{0 < |\bm{k}| \le R} |\bm{k}|^{-2} \right)\left(\sum_{0 < |\bm{k}| \le R} |\bm{k}|^{2}\left|\hat{f}_{\bm{k}}\right|^2 \right)\leq \frac{RL^2}{4\pi^2} \left( \sum_{0 < |\bm{k}| \le R} \mu^2|\bm{k}|^{2}\left|\hat{f}_{\bm{k}}\right|^2 \right) \leq \frac{R}{4\pi^2\del{ }\ins{\thinspace }L} \left\Vert \nabla_N \bm f \right\Vert^2
\end{aligned}
\]
and
\[
\left(\sum_{|\bm{k}| > R} \left|\hat{f}_{\bm{k}}\right| \right)^2 \leq \left(\sum_{|\bm{k}| > R} |\bm{k}|^{-4} \right)\left(\sum_{|\bm{k}| > R} |\bm{k}|^{4}\left|\hat{f}_{\bm{k}}\right|^2 \right) \leq \frac{L}{16\pi^4 R} \left\Vert \Delta_N \bm f \right\Vert^2,
\]
respectively.
Next, choosing $R = (L \left\Vert \Delta_N \bm f \right\Vert_2) / (2\pi \left\Vert \nabla_N \bm f \right\Vert_2)$ and substituting the above two inequalities into \eqref{eqn:infexpend} yields
\[
\left\Vert \bm f \right\Vert_\infty^2 \leq 3 |\Omega|^{-2} |\bar{\bm f}|^2  +  \frac{3}{8\pi^3} \left\Vert \nabla_N \bm f \right\Vert_2 \left\Vert \Delta_N \bm f \right\Vert_2.
\]
Thus, the first inequality \eqref{eqn:3dSob} is established. With the help of Fourier expansion, the second inequality \eqref{eqn:poincareIneq} can be directly verified.

{\opecolor{red}
To prove the third inequality \eqref{eqn:nonlineEst1}, by using \cite[Lemma 1]{gottlieb2012stability}, one has
\[
 \left\Vert \Delta_N (\bm f \bm g) \right\Vert_2^2  =  C \left\Vert \Delta I_N ( \bm f \bm g) \right\Vert_{L^2}^2 = C\left\Vert \Delta (f g) \right\Vert_{L^2}^2,
\]
where $I_N$ denotes the interpolation operator mapping from trigonometric space to $L^2$ space.
By the chain rule for the Laplacian operator and the triangle inequality, we have
\[
\left\Vert \Delta (f g) \right\Vert_{L^2}^2 \leq 3(\left\Vert g \Delta f \right\Vert_{L^2}^2 + \left\Vert f \Delta g \right\Vert_{L^2}^2  + 2 \left\Vert  \nabla f  \cdot \nabla g \right\Vert_{L^2}^2).
\]
For the first two leading terms on the right-hand side, we obtain
\[
    \left\Vert g \Delta f \right\Vert_{L^2}^2 \leq \left \Vert g \right\Vert_{L^\infty}^2 \left\Vert  \Delta f \right\Vert_{L^2}^2,\quad  \left\Vert f \Delta g \right\Vert_{L^2}^2 \leq \left \Vert f \right\Vert_{L^\infty}^2 \left\Vert \Delta g \right\Vert_{L^2}^2.
\]
For the cross term, we first apply the H\"older inequality and the Sobolev interpolation inequality, and subsequently employ Young's inequality to obtain
\begin{align*}
\left\Vert  \nabla f  \cdot \nabla g \right\Vert_{L^2}^2 \leq & \left\Vert  \nabla f \right\Vert_{L^4}^2 \left\Vert  \nabla g \right\Vert_{L^4}^2  \leq \left \Vert  f \right\Vert_{L^\infty} \left\Vert \Delta f \right\Vert_{L^2} \left \Vert g \right\Vert_{L^\infty}  \left\Vert  \Delta  g \right\Vert_{L^2} \\
\leq & \frac{1}{2} \left \Vert  g \right\Vert_{L^\infty}^2 \left\Vert \Delta f \right\Vert_{L^2}^2 + \frac{1}{2}  \left \Vert f \right\Vert_{L^\infty}^2  \left\Vert  \Delta  g \right\Vert_{L^2}^2.
\end{align*}
Combining all estimates, we get
\[
\begin{aligned}
     \left\Vert \Delta_N (\bm f \bm g) \right\Vert_2^2  & \leq  C\left\Vert \Delta (f g) \right\Vert_{L^2}^2 \leq C( \left \Vert  g \right\Vert_{L^\infty}^2 \left\Vert \Delta f \right\Vert_{L^2}^2 + \left \Vert f \right\Vert_{L^\infty}^2  \left\Vert  \Delta  g \right\Vert_{L^2}^2)  \\
     &\leq C( \left \Vert \bm g \right\Vert_{\infty}^2 \left\Vert \Delta_N \bm  f \right\Vert_{2}^2 + \left \Vert \bm  f \right\Vert_{\infty}^2  \left\Vert \Delta_N \bm g \right\Vert_{2}^2),
\end{aligned}
\]
where the last inequality holds via the interpolation equivalence. Thus, taking the square root of both sides yields \eqref{eqn:nonlineEst1}.  The inequality \eqref{eqn:nonlineEst2} can be derived in a similar manner. Moreover, \eqref{eqn:nonlineEst3} follows immediately by recursively using \eqref{eqn:nonlineEst1} and \eqref{eqn:nonlineEst2}. This completes the proof.
}
\end{appgroup}

\bibliographystyle{cas-model1-sc}

\begin{thebibliography}{34}
\expandafter\ifx\csname natexlab\endcsname\relax\def\natexlab#1{#1}\fi
\ifx\bibinfo\relax\providecommand{\bibinfo}[2]{#2}\fi
\renewcommand{\newblock}{}
\makeatletter
\ifx\xfnm\@undefined\def\xfnm[#1]{#1}\fi
\ifx\xsnm\@undefined\def\xsnm[#1]{#1}\fi
\ifx\plxcitation\@undefined\def\plxcitation#1#2#3#4#5#6{#6}\fi
\ifx\endplxcitation\@undefined\def\endplxcitation{}\fi
\ifx\PrintOrdinal\@undefined\def\PrintOrdinal#1{#1}\fi
\makeatother
\bibitem[{Cahn and Hilliard(1958)}]{cahn1958free}
\plxcitation{}{cahn1958free}{sbref0001}{}{article}{%
\bibinfo{author}{\xfnm[J.~W.] \xsnm[Cahn]}, \bibinfo{author}{\xfnm[J.~E.]
  \xsnm[Hilliard]},
\newblock \bibinfo{title}{{{Free energy of a nonuniform system. I. \del{I}\ins{i}nterfacial
  free energy}}},
\newblock \bibinfo{journal}{\del{The Journal of Chemical Physics}\ins{J. Chem. Phys.}}
  \bibinfo{volume}{28} (\bibinfo{number}{2}) (\bibinfo{year}{1958})
  \bibinfo{pages}{258--267}.
}\endplxcitation
\bibitem[{Eyre(1998)}]{eyre1998unconditionally-ch}
\plxcitation{}{eyre1998unconditionally-ch}{sbref0002}{}{inproceedings}{%
\bibinfo{author}{\xfnm[D.~J.] \xsnm[Eyre]},
\newblock \bibinfo{title}{{{Unconditionally gradient stable time marching the
  Cahn-Hilliard equation}}},
\newblock \bibinfo{booktitle}{MRS Online Proceedings Library}
  \bibinfo{volume}{529} (\bibinfo{year}{1998}) \bibinfo{pages}{39}.
}\endplxcitation
\bibitem[{Wise(2010)}]{wise2010unconditionally}
\plxcitation{}{wise2010unconditionally}{sbref0003}{}{article}{%
\bibinfo{author}{\xfnm[S.~M.] \xsnm[Wise]},
\newblock \bibinfo{title}{{{Unconditionally stable finite difference, nonlinear
  multigrid simulation of the Cahn-Hilliard-Hele-Shaw system of equations}}},
\newblock \bibinfo{journal}{\del{Journal of Scientific Computing}\ins{J. Sci. Comput.}}
  \bibinfo{volume}{44} (\bibinfo{number}{1}) (\bibinfo{year}{2010})
  \bibinfo{pages}{38--68}.
}\endplxcitation
\bibitem[{Guo {\rm et~al.}(2024)Guo, Wang, Wise, and Zhang}]{GuoWanWise24}
\plxcitation{}{GuoWanWise24}{sbref0004}{}{article}{%
\bibinfo{author}{\xfnm[Y.] \xsnm[Guo]}, \bibinfo{author}{\xfnm[C.]
  \xsnm[Wang]}, \bibinfo{author}{\xfnm[S.] \xsnm[Wise]},
  \bibinfo{author}{\xfnm[Z.] \xsnm[Zhang]},
\newblock \bibinfo{title}{{{Convergence analysis of a positivity-preserving
  numerical scheme for the Cahn-Hilliard-Stokes system with Flory-Huggins
  energy potential}}},
\newblock \bibinfo{journal}{\del{Mathematics of Computation}\ins{Math. Comput.}} \bibinfo{volume}{93}
  (\bibinfo{number}{349}) (\bibinfo{year}{2024}) \bibinfo{pages}{2185--2214}.
}\endplxcitation
\bibitem[{Guo {\rm et~al.}(2026)Guo, Wang, Wise, and Zhang}]{GuoWanWise26}
\plxcitation{}{GuoWanWise26}{sbref0005}{}{article}{%
\bibinfo{author}{\xfnm[Y.] \xsnm[Guo]}, \bibinfo{author}{\xfnm[C.]
  \xsnm[Wang]}, \bibinfo{author}{\xfnm[S.~M.] \xsnm[Wise]},
  \bibinfo{author}{\xfnm[Z.] \xsnm[Zhang]},
\newblock \bibinfo{title}{{{A uniquely solvable and positivity-preserving finite
  difference scheme for the Flory--Huggins--Cahn--Hilliard equation with
  dynamical boundary condition}}},
\newblock \bibinfo{journal}{\del{Journal of Computational and Applied Mathematics}\ins{J. Comput. Appl. Math.}}
  \bibinfo{volume}{472} (\bibinfo{year}{2026}) \bibinfo{pages}{116789}.
}\endplxcitation
\bibitem[{Li {\rm et~al.}(2024)Li, Jing, Liu, Wang, and Chen}]{LiLiuWan24}
\plxcitation{}{LiLiuWan24}{sbref0006}{}{article}{%
\bibinfo{author}{\xfnm[Y.] \xsnm[Li]}, \bibinfo{author}{\xfnm[J.] \xsnm[Jing]},
  \bibinfo{author}{\xfnm[Q.] \xsnm[Liu]}, \bibinfo{author}{\xfnm[C.]
  \xsnm[Wang]}, \bibinfo{author}{\xfnm[W.] \xsnm[Chen]},
\newblock \bibinfo{title}{{{A third order positivity-preserving, energy stable
  numerical scheme for the Cahn-Hilliard equation with logarithmic potential}}},
\newblock \bibinfo{journal}{Sci. Sin. Math.} \bibinfo{volume}{55} (\bibinfo{number}{5}) (\bibinfo{year}{2025}) \bibinfo{pages}{991--1020} (in Chinese).
}\endplxcitation
\bibitem[{Sun {\rm et~al.}(2023)Sun, Zhang, Qian, and Song}]{SunZhanQian}
\plxcitation{}{SunZhanQian}{sbref0007}{}{article}{%
\bibinfo{author}{\xfnm[J.] \xsnm[Sun]}, \bibinfo{author}{\xfnm[H.]
  \xsnm[Zhang]}, \bibinfo{author}{\xfnm[X.] \xsnm[Qian]},
  \bibinfo{author}{\xfnm[S.] \xsnm[Song]},
\newblock \bibinfo{title}{{{A family of structure-preserving exponential time
  differencing Runge--Kutta schemes for the viscous Cahn--Hilliard equation}}},
\newblock \bibinfo{journal}{\del{Journal of Computational Physics}\ins{J. Comput. Phys.}}
  \bibinfo{volume}{492} (\bibinfo{year}{2023}) \bibinfo{pages}{112414}.
}\endplxcitation
\bibitem[{Li {\rm et~al.}(2022)Li, Quan, and Tang}]{LiQuanTang}
\plxcitation{}{LiQuanTang}{sbref0008}{}{article}{%
\bibinfo{author}{\xfnm[D.] \xsnm[Li]}, \bibinfo{author}{\xfnm[C.] \xsnm[Quan]},
  \bibinfo{author}{\xfnm[T.] \xsnm[Tang]},
\newblock \bibinfo{title}{{{Stability and convergence analysis for the
  implicit-explicit method to the Cahn-Hilliard equation}}},
\newblock \bibinfo{journal}{\del{Mathematics of Computation}\ins{Math. Comput.}} \bibinfo{volume}{91}
  (\bibinfo{number}{334}) (\bibinfo{year}{2022}) \bibinfo{pages}{785--809}.
}\endplxcitation
\bibitem[{Li and Quan(2022)}]{LiQuan}
\plxcitation{}{LiQuan}{sbref0009}{}{article}{%
\bibinfo{author}{\xfnm[D.] \xsnm[Li]}, \bibinfo{author}{\xfnm[C.] \xsnm[Quan]},
\newblock \bibinfo{title}{{{The operator-splitting method for Cahn-Hilliard is
  stable}}},
\newblock \bibinfo{journal}{\del{Journal of Scientific Computing}\ins{J. Sci. Comput.}}
  \bibinfo{volume}{90} (\bibinfo{number}{1}) (\bibinfo{year}{2022})
  \bibinfo{pages}{62}.
}\endplxcitation
\bibitem[{Furihata(2001)}]{furihata2001stable}
\plxcitation{}{furihata2001stable}{sbref0010}{}{article}{%
\bibinfo{author}{\xfnm[D.] \xsnm[Furihata]},
\newblock \bibinfo{title}{{{A stable and conservative finite difference scheme
  for the Cahn-Hilliard equation}}},
\newblock \bibinfo{journal}{\del{Numerische Mathematik}\ins{Numer. Math.}} \bibinfo{volume}{87}
  (\bibinfo{number}{4}) (\bibinfo{year}{2001}) \bibinfo{pages}{675--699}.
}\endplxcitation
\bibitem[{Zhang {\rm et~al.}(2024)Zhang, Liu, Qian, and Song}]{ZhanLiu24}
\plxcitation{}{ZhanLiu24}{sbref0011}{}{article}{%
\bibinfo{author}{\xfnm[H.] \xsnm[Zhang]}, \bibinfo{author}{\xfnm[L.]
  \xsnm[Liu]}, \bibinfo{author}{\xfnm[X.] \xsnm[Qian]},
  \bibinfo{author}{\xfnm[S.] \xsnm[Song]},
\newblock \bibinfo{title}{{{Large time-stepping, delay-free, and
  invariant-set-preserving integrators for the viscous Cahn--Hilliard--Oono
  equation}}},
\newblock \bibinfo{journal}{\del{Journal of Computational Physics}\ins{J. Comput. Phys.}}
  \bibinfo{volume}{499} (\bibinfo{year}{2024}) \bibinfo{pages}{112708}.
}\endplxcitation
\bibitem[{Chen {\rm et~al.}(2024)Chen, Jing, Liu, Wang, and Wang}]{chen2024second}
\plxcitation{}{chen2024second}{sbref0012}{}{article}{%
\bibinfo{author}{\xfnm[W.] \xsnm[Chen]}, \bibinfo{author}{\xfnm[J.]
  \xsnm[Jing]}, \bibinfo{author}{\xfnm[Q.] \xsnm[Liu]},
  \bibinfo{author}{\xfnm[C.] \xsnm[Wang]}, \bibinfo{author}{\xfnm[X.]
  \xsnm[Wang]},
\newblock \bibinfo{title}{{{A second order numerical scheme of the Cahn-Hilliard-Navier-Stokes system with Flory-Huggins potential}}},
\newblock \bibinfo{journal}{\del{Communications in computational physics}\ins{Commun. Comput. Phys.}}
\bibinfo{volume}{35} (\bibinfo{number}{3}) (\bibinfo{year}{2024}) \bibinfo{pages}{633--661}.}\endplxcitation
\bibitem[{Cheng {\rm et~al.}(2022)Cheng, Wang, Wise, and Wu}]{cheng2022third}
\plxcitation{}{cheng2022third}{sbref0013}{}{article}{%
\bibinfo{author}{\xfnm[K.] \xsnm[Cheng]}, \bibinfo{author}{\xfnm[C.]
  \xsnm[Wang]}, \bibinfo{author}{\xfnm[S.~M.] \xsnm[Wise]},
  \bibinfo{author}{\xfnm[Y.] \xsnm[Wu]},
\newblock \bibinfo{title}{{{A third order accurate in time, {BDF}-type energy
  stable scheme for the Cahn-Hilliard equation}}},
\newblock \bibinfo{journal}{\del{Numerical Mathematics: Theory, Methods and
  Applications}\ins{Numer. Math. Theory Methods Appl.}} \bibinfo{volume}{15} (\bibinfo{number}{2})
  (\bibinfo{year}{2022}) \bibinfo{pages}{279--303}.
}\endplxcitation
\bibitem[{Feng and Prohl(2003)}]{feng2003numerical}
\plxcitation{}{feng2003numerical}{sbref0014}{}{article}{%
\bibinfo{author}{\xfnm[X.] \xsnm[Feng]}, \bibinfo{author}{\xfnm[A.]
  \xsnm[Prohl]},
\newblock \bibinfo{title}{{{Numerical analysis of the Allen-Cahn equation and
  approximation for mean curvature flows}}},
\newblock \bibinfo{journal}{\del{Numerische Mathematik}\ins{Numer. Math.}} \bibinfo{volume}{94}
  (\bibinfo{number}{1}) (\bibinfo{year}{2003}) \bibinfo{pages}{33--65}.
}\endplxcitation
\bibitem[{Meng {\rm et~al.}(2025)Meng, Bao, and Zhang}]{xiangjun2025second}
\plxcitation{}{xiangjun2025second}{sbref0015}{}{article}{%
\bibinfo{author}{\xfnm[X.] \xsnm[Meng]}, \bibinfo{author}{\xfnm[X.]
  \xsnm[Bao]}, \bibinfo{author}{\xfnm[Z.] \xsnm[Zhang]},
\newblock \bibinfo{title}{{{Second-order linear stabilized semi-implicit
  Crank-Nicolson scheme for the Cahn-Hilliard model with dynamic boundary
  conditions}}},
\newblock \bibinfo{journal}{\del{Communications in Computational Physics}\ins{Commun. Comput. Phys.}}
  \bibinfo{volume}{37} (\bibinfo{number}{1}) (\bibinfo{year}{2025})
  \bibinfo{pages}{137--170}.
}\endplxcitation
\bibitem[{Li {\rm et~al.}(2023)Li, Qiao, and Wang}]{Li2023Stabilization}
\plxcitation{}{Li2023Stabilization}{sbref0016}{}{article}{%
\bibinfo{author}{\xfnm[X.] \xsnm[Li]}, \bibinfo{author}{\xfnm[Z.] \xsnm[Qiao]},
  \bibinfo{author}{\xfnm[C.] \xsnm[Wang]},
\newblock \bibinfo{title}{{{Stabilization parameter analysis of a second-order
  linear numerical scheme for the nonlocal Cahn--Hilliard equation}}},
\newblock \bibinfo{journal}{\del{IMA Journal of Numerical Analysis}\ins{IMA J. Numer. Anal.}}
  \bibinfo{volume}{43} (\bibinfo{year}{2023}) \bibinfo{pages}{1089--1114}.
}\endplxcitation
\bibitem[{Fu and Yang(2022)}]{fu2022energy}
\plxcitation{}{fu2022energy}{sbref0017}{}{article}{%
\bibinfo{author}{\xfnm[Z.] \xsnm[Fu]}, \bibinfo{author}{\xfnm[J.] \xsnm[Yang]},
\newblock \bibinfo{title}{{{Energy-decreasing exponential time differencing
  {R}unge--{K}utta methods for phase-field models}}},
\newblock \bibinfo{journal}{\del{Journal of Computational Physics}\ins{J. Comput. Phys.}}
  \bibinfo{volume}{454} (\bibinfo{year}{2022}) \bibinfo{pages}{110943}.
}\endplxcitation
\bibitem[{Fu {\rm et~al.}(2024)Fu, Tang, and Yang}]{fu2022unconditionally}
\plxcitation{}{fu2022unconditionally}{sbref0018}{}{article}{%
\bibinfo{author}{\xfnm[Z.] \xsnm[Fu]}, \bibinfo{author}{\xfnm[T.] \xsnm[Tang]},
  \bibinfo{author}{\xfnm[J.] \xsnm[Yang]},
\newblock \bibinfo{title}{{{Unconditionally energy decreasing high-order
  \del{I}\ins{i}mplicit-explicit Runge-Kutta methods for phase-field models with the
  Lipschitz nonlinearity}}},
\newblock \bibinfo{journal}{\del{Mathematics of Computation}\ins{Math. Comput.}} \bibinfo{volume}{93}
  (\bibinfo{year}{2024}) \bibinfo{pages}{2745--2767}.
}\endplxcitation
\bibitem[{Li and Qiao(2017{\natexlab{a}})}]{li2017second}
\plxcitation{}{li2017second}{sbref0019}{}{article}{%
\bibinfo{author}{\xfnm[D.] \xsnm[Li]}, \bibinfo{author}{\xfnm[Z.] \xsnm[Qiao]},
\newblock \bibinfo{title}{{{On second order semi-implicit Fourier spectral
  methods for 2D Cahn--Hilliard equations}}},
\newblock \bibinfo{journal}{\del{Journal of scientific computing}\ins{J. Sci. Comput.}}
  \bibinfo{volume}{70} (\bibinfo{year}{2017}{\natexlab{a}})
  \bibinfo{pages}{301--341}.
}\endplxcitation
\bibitem[{Li and Qiao(2017{\natexlab{b}})}]{li2017stabilization}
\plxcitation{}{li2017stabilization}{sbref0020}{}{article}{%
\bibinfo{author}{\xfnm[D.] \xsnm[Li]}, \bibinfo{author}{\xfnm[Z.] \xsnm[Qiao]},
\newblock \bibinfo{title}{{{On the stabilization size of semi-implicit
  Fourier-spectral methods for 3D Cahn--Hilliard equations}}},
\newblock \bibinfo{journal}{\del{Communications in Mathematical Sciences}\ins{Commun. Math. Sci.}}
  \bibinfo{volume}{15} (\bibinfo{number}{6})
  (\bibinfo{year}{2017}{\natexlab{b}}) \bibinfo{pages}{1489--1506}.
}\endplxcitation
\bibitem[{Li and Qiao(2024)}]{LiQiao24}
\plxcitation{}{LiQiao24}{sbref0021}{}{article}{%
\bibinfo{author}{\xfnm[X.] \xsnm[Li]}, \bibinfo{author}{\xfnm[Z.] \xsnm[Qiao]},
\newblock \bibinfo{title}{{A second-order, linear, $L^\infty$-convergent, and
  energy stable scheme for the phase field crystal equation}},
\newblock \bibinfo{journal}{\del{SIAM Journal on Scientific Computing}\ins{SIAM J. Sci. Comput.}}
  \bibinfo{volume}{46} (\bibinfo{number}{1}) (\bibinfo{year}{2024})
  \bibinfo{pages}{A429--A451}.
}\endplxcitation
\bibitem[{Gottlieb {\rm et~al.}(2012)Gottlieb, Tone, Wang, Wang, and
  Wirosoetisno}]{gottlieb2012long}
\plxcitation{}{gottlieb2012long}{sbref0022}{}{article}{%
\bibinfo{author}{\xfnm[S.] \xsnm[Gottlieb]}, \bibinfo{author}{\xfnm[F.]
  \xsnm[Tone]}, \bibinfo{author}{\xfnm[C.] \xsnm[Wang]},
  \bibinfo{author}{\xfnm[X.] \xsnm[Wang]}, \bibinfo{author}{\xfnm[D.]
  \xsnm[Wirosoetisno]},
\newblock \bibinfo{title}{{{Long time stability of a classical efficient scheme
  for two-dimensional Navier--Stokes equations}}},
\newblock \bibinfo{journal}{\del{SIAM Journal on Numerical Analysis}\ins{SIAM J. Numer. Anal.}}
  \bibinfo{volume}{50} (\bibinfo{number}{1}) (\bibinfo{year}{2012})
  \bibinfo{pages}{126--150}.
}\endplxcitation
\bibitem[{Li {\rm et~al.}(2025)Li, Qiao, Wang, and Zheng}]{li2025global}
\plxcitation{}{li2025global}{sbref0023}{}{article}{%
\bibinfo{author}{\xfnm[X.] \xsnm[Li]}, \bibinfo{author}{\xfnm[Z.] \xsnm[Qiao]},
  \bibinfo{author}{\xfnm[C.] \xsnm[Wang]}, \bibinfo{author}{\xfnm[N.]
  \xsnm[Zheng]},
\newblock \bibinfo{title}{{{Global-in-time energy stability analysis for a
  second-order accurate exponential time differencing Runge--Kutta scheme for
  the phase field crystal equation}}},
\newblock \bibinfo{journal}{\del{Mathematics of Computation}\ins{Math. Comput.}}
  \bibinfo{volume}{95} (\bibinfo{number}{358}) (\bibinfo{year}{2026}) \bibinfo{pages}{803--831}.}\endplxcitation
\bibitem[{Zhang {\rm et~al.}(2024)Zhang, Wang, and Teng}]{zhang2024second}
\plxcitation{}{zhang2024second}{sbref0024}{}{article}{%
\bibinfo{author}{\xfnm[H.] \xsnm[Zhang]}, \bibinfo{author}{\xfnm[H.]
  \xsnm[Wang]}, \bibinfo{author}{\xfnm[X.] \xsnm[Teng]},
\newblock \bibinfo{title}{{{A second-order, global-in-time energy stable
  implicit-explicit Runge--Kutta scheme for the phase field crystal equation}}},
\newblock \bibinfo{journal}{\del{SIAM Journal on Numerical Analysis}\ins{SIAM J. Numer. Anal.}}
  \bibinfo{volume}{62} (\bibinfo{number}{6}) (\bibinfo{year}{2024})
  \bibinfo{pages}{2667--2697}.
}\endplxcitation
\bibitem[{Cheng {\rm et~al.}(2016)Cheng, Wang, Wise, and Yue}]{cheng2016second}
\plxcitation{}{cheng2016second}{sbref0025}{}{article}{%
\bibinfo{author}{\xfnm[K.] \xsnm[Cheng]}, \bibinfo{author}{\xfnm[C.]
  \xsnm[Wang]}, \bibinfo{author}{\xfnm[S.~M.] \xsnm[Wise]},
  \bibinfo{author}{\xfnm[X.] \xsnm[Yue]},
\newblock \bibinfo{title}{{{A second-order, weakly energy-stable pseudo-spectral
  scheme for the Cahn--Hilliard equation and its solution by the homogeneous
  linear iteration method}}},
\newblock \bibinfo{journal}{\del{Journal of Scientific Computing}\ins{J. Sci. Comput.}}
  \bibinfo{volume}{69} (\bibinfo{number}{3}) (\bibinfo{year}{2016})
  \bibinfo{pages}{1083--1114}.
}\endplxcitation
\bibitem[{Cheng(2019)}]{cheng2019energy}
\plxcitation{}{cheng2019energy}{sbref0026}{}{article}{%
\bibinfo{author}{\xfnm[K.] \xsnm[Cheng]},
\newblock \bibinfo{title}{{An \del{E}\ins{e}nergy \del{S}\ins{s}table {BDF2} {F}ourier \del{P}\ins{p}seudo-spectral
  \del{N}\ins{n}umerical \del{S}\ins{s}cheme for the \del{S}\ins{s}quare \del{P}\ins{p}hase \del{F}\ins{f}ield \del{C}\ins{c}rystal \del{E}\ins{e}quation}},
\newblock \bibinfo{journal}{\del{Communications in Computational Physics}\ins{Commun. Comput. Phys.}}
  \bibinfo{volume}{26} (\bibinfo{number}{5}) (\bibinfo{year}{2019})
  \bibinfo{pages}{1335--1364}.
}\endplxcitation
\bibitem[{Higham(2008)}]{higham2008functions}
\plxcitation{}{higham2008functions}{sbref0027}{}{book}{%
\bibinfo{author}{\xfnm[N.~J.] \xsnm[Higham]}, \bibinfo{title}{{Functions of
  Matrices: Theory and Computation}}, \bibinfo{publisher}{SIAM},
  \bibinfo{year}{2008}.
}\endplxcitation
\bibitem[{Li(2019)}]{li2019convergence}
\plxcitation{}{li2019convergence}{sbref0028}{}{article}{%
\bibinfo{author}{\xfnm[X.] \xsnm[Li]}, \bibinfo{author}{\xfnm[L.] \xsnm[Ju]}, \bibinfo{author}{\xfnm[X.] \xsnm[Meng]},
\newblock \bibinfo{title}{{{Convergence analysis of exponential time
  differencing schemes for the Cahn-Hilliard equation}}},
\newblock \bibinfo{journal}{\del{Communications in Computational Physics}\ins{Commun. Comput. Phys.}}
  \bibinfo{volume}{26} (\bibinfo{number}{5}) (\bibinfo{year}{2019}) \bibinfo{pages}{1510--1529}.}\endplxcitation
\bibitem[{Hochbruck and Ostermann(2005)}]{hochbruck2005explicit}
\plxcitation{}{hochbruck2005explicit}{sbref0029}{}{article}{%
\bibinfo{author}{\xfnm[M.] \xsnm[Hochbruck]}, \bibinfo{author}{\xfnm[A.]
  \xsnm[Ostermann]},
\newblock \bibinfo{title}{{{Explicit exponential Runge--Kutta methods for
  semilinear parabolic problems}}},
\newblock \bibinfo{journal}{\del{SIAM Journal on Numerical Analysis}\ins{SIAM J. Numer. Anal.}}
  \bibinfo{volume}{43} (\bibinfo{number}{3}) (\bibinfo{year}{2005})
  \bibinfo{pages}{1069--1090}.
}\endplxcitation
\bibitem[{Elliott and Songmu(1986)}]{elliott1986cahn}
\plxcitation{}{elliott1986cahn}{sbref0030}{}{article}{%
\bibinfo{author}{\xfnm[C.~M.] \xsnm[Elliott]}, \bibinfo{author}{\xfnm[Z.]
  \xsnm[Songmu]},
\newblock \bibinfo{title}{{{On the Cahn-Hilliard equation}}},
\newblock \bibinfo{journal}{\del{Archive for Rational Mechanics and Analysis}\ins{Arch. Ration. Mech. Anal.}}
  \bibinfo{volume}{96} (\bibinfo{number}{4}) (\bibinfo{year}{1986})
  \bibinfo{pages}{339--357}.
}\endplxcitation
\bibitem[{Novick-Cohen(2008)}]{novick2008cahn}
\plxcitation{}{novick2008cahn}{sbref0031}{}{article}{%
\bibinfo{author}{\xfnm[A.] \xsnm[Novick-Cohen]},
\newblock \bibinfo{title}{{{The Cahn--Hilliard equation}}},
\newblock \bibinfo{journal}{\del{Handbook of differential equations: evolutionary
  equations}\ins{Handb. Differ. Eqs. Evol. Eqs.}} \bibinfo{volume}{4} (\bibinfo{year}{2008})
  \bibinfo{pages}{201--228}.
}\endplxcitation
\bibitem[{Temam(2012)}]{temam2012infinite}
\plxcitation{}{temam2012infinite}{sbref0032}{}{book}{%
\bibinfo{author}{\xfnm[R.] \xsnm[Temam]}, \bibinfo{title}{Infinite-Dimensional
  Dynamical Systems in Mechanics and Physics}, \bibinfo{publisher}{Springer
  Science \& Business Media}, \bibinfo{year}{2012}.
}\endplxcitation
\bibitem[{Cheng and Wang(2016)}]{cheng2016long}
\plxcitation{}{cheng2016long}{sbref0033}{}{article}{%
\bibinfo{author}{\xfnm[K.] \xsnm[Cheng]}, \bibinfo{author}{\xfnm[C.]
  \xsnm[Wang]},
\newblock \bibinfo{title}{{{Long time stability of high order multistep
  numerical schemes for two-dimensional incompressible Navier--Stokes
  equations}}},
\newblock \bibinfo{journal}{\del{SIAM Journal on Numerical Analysis}\ins{SIAM J. Numer. Anal.}}
  \bibinfo{volume}{54} (\bibinfo{number}{5}) (\bibinfo{year}{2016})
  \bibinfo{pages}{3123--3144}.
}\endplxcitation
\bibitem[{Gottlieb and Wang(2012)}]{gottlieb2012stability}
\plxcitation{}{gottlieb2012stability}{sbref0034}{}{article}{%
\bibinfo{author}{\xfnm[S.] \xsnm[Gottlieb]}, \bibinfo{author}{\xfnm[C.]
  \xsnm[Wang]},
\newblock \bibinfo{title}{{Stability and \del{C}\ins{c}onvergence \del{A}\ins{a}nalysis of \del{F}\ins{f}ully \del{D}\ins{d}iscrete
  {F}ourier \del{C}\ins{c}ollocation \del{S}\ins{s}pectral \del{M}\ins{m}ethod for 3-{D} \del{V}\ins{v}iscous {Burgers'} \del{E}\ins{e}quation}},
\newblock \bibinfo{journal}{\del{Journal of Scientific Computing}\ins{J. Sci. Comput.}}
  \bibinfo{volume}{53} (\bibinfo{number}{1}) (\bibinfo{year}{2012})
  \bibinfo{pages}{102--128}.
  \bibinfo{doi}{10.1007/s10915-012-9621-8}}\endplxcitation
\end{thebibliography}

\end{document}